\documentclass[aos,MSNbibl,seceqn,dvips]{arximspdf}
\usepackage{graphicx}
\doi{10.1214/14-AOS1213} 
\volume{42}
\issue{3}
\pubyear{2014}
\firstpage{1029}
\lastpage{1069}

\makeatletter
\newcommand{\rrvert}{\vert}
\newcommand{\llvert}{\vert}
\renewcommand{\backslash}{\setminus}
\newcommand{\eqref}[1]{(\ref{#1})}

\newtheorem{theo}{Theorem}
\newproclaim{rema}[theo]{Remark}
\newtheorem{lem}[theo]{Lemma}
\newproclaim{Assumption}{Assumption}
\newproclaim{Example}{Example}

\newcommand{\R}{\mathbb{R}}
\newcommand{\E}{\mathbb{E}}
\newcommand{\PP}{\mathbb{P}}
\newcommand{\al}{\alpha}
\newcommand{\Ga}{\Gamma}
\newcommand{\ka}{\kappa}
\newcommand{\vpi}{\varpi}
\newcommand{\si}{\sigma}
\newcommand{\Te}{\Theta}
\newcommand{\be}{\beta}
\newcommand{\ep}{\varepsilon}
\newcommand{\De}{\Delta}
\newcommand{\Om}{\Omega}
\newcommand{\ze}{\zeta}
\newcommand{\f}{\mathcal{F}}
\newcommand{\g}{\mathcal{G}}
\newcommand{\laa}{\mathcal{L}}
\newcommand{\ua}{\mathcal{U}}
\newcommand{\wC}{\widehat{C}}
\newcommand{\wY}{\widehat{Y}}
\newcommand{\wU}{\widehat{U}}
\newcommand{\wde}{\widehat{\delta}}
\newcommand{\WE}{\widetilde{\mathbb{E}}}
\newcommand{\WZ}{\widetilde{Z}}
\newcommand{\BF}{\overline{F}{}}
\newcommand{\BR}{\overline{R}{}}
\newcommand{\BU}{\overline{U}{}}
\newcommand{\BV}{\overline{V}{}}
\newcommand{\BZ}{\overline{Z}{}}
\newcommand{\BX}{\overline{X}{}}
\newcommand{\Bte}{\overline{\theta}{}}
\newcommand{\Beta}{\overline{\eta}{}}
\newcommand{\Brho}{\overline{\rho}{}}
\newcommand{\rdn}{\sqrt{\De_n}}
\newcommand{\um}{p\!\!\!\!{}_{_=}\!\,}
\newcommand{\un}{q\!\!\!{}_{_=}\!\,}
\newcommand{\toop}{\stackrel{\PP}{\longrightarrow}}
\newcommand{\tolls}{\stackrel{\laa- s}{\Longrightarrow}}
\newcommand{\tols}{\stackrel{\laa-s}{\longrightarrow}}
\newcommand{\toucp}{\stackrel{{\mathrm{u.c.p.}}}{\Longrightarrow}}
\newcommand{\sign}{\operatorname{sign}}
\makeatother

\begin{document}
\begin{frontmatter}

\title{Efficient estimation of integrated volatility in presence of
infinite variation jumps}
\runtitle{Efficient estimation of integrated volatility}

\begin{aug}
\author[a]{\fnms{Jean}~\snm{Jacod}\ead[label=e1]{jean.jacod@upmc.fr}}
\and
\author[b]{\fnms{Viktor}~\snm{Todorov}\corref{}\thanksref{t1}\ead[label=e2]{v-todorov@northwestern.edu}}
\runauthor{J. Jacod and V. Todorov}
\affiliation{UPMC (Universit\'e Paris-6) and Northwestern University}
\address[a]{Institut de Math\'ematiques de Jussieu\\
CNRS-UMR 7586\\
Universit\'e Pierre et Marie Curie--P6\\
4 Place Jussieu\\
75252 Paris-Cedex 05\\
France\\
\printead{e1}} 
\address[b]{Department of Finance\\
Northwestern University\\
Evanston, Illinois 60208-2001\\
USA\\
\printead{e2}}
\end{aug}
\thankstext{t1}{Supported in part by NSF Grant SES-0957330.}

\received{\smonth{6} \syear{2013}}
\revised{\smonth{12} \syear{2013}}

%
\begin{abstract}
We propose new nonparametric estimators of the integrated volatility of
an It\^o semimartingale observed at discrete times on a fixed time interval
with mesh of the observation grid shrinking to zero. The proposed estimators
achieve the optimal rate and variance of estimating integrated
volatility even
in the presence of infinite variation jumps when the latter are stochastic
integrals with respect to locally ``stable'' L\'evy processes, that is,
processes
whose L\'evy measure around zero behaves like that of a stable process.
On a
first step, we estimate locally volatility from the empirical characteristic
function of the increments of the process over blocks of shrinking
length and
then we sum these estimates to form initial estimators of the integrated
volatility. The estimators contain bias when jumps of infinite
variation are
present, and on a second step we estimate and remove this bias by using
integrated volatility estimators formed from the empirical characteristic
function of the high-frequency increments for different values of its
argument. The second step debiased estimators achieve efficiency and we
derive a feasible central limit theorem for them.
\end{abstract}

%

\end{frontmatter}

\section{Introduction}\label{sec-Intro}

In this paper, we consider the problem of estimating the continuous
part of
the quadratic variation (henceforth referred to as integrated
volatility) of
a discretely-observed one-dimensional It\^o semimartingale over a finite
interval with mesh of the observation grid going to zero in the case
when the
observed process can contain jumps of infinite variation. Separating jumps
from diffusive volatility is of central interest in finance due to the
distinct role played by diffusive volatility and jumps in financial decision
making, which is also reflected in the distinct risk premium demanded by
investors for each of them; see, for example, \cite{BT}. Until now,
this problem has
been well studied when jumps are of finite variation; see, for example,
\cite{BS,BSW,M1,M2,JP}. However, empirical results in
\cite{AJ}
suggest that for some financial data sets jumps can be of infinite variation.
This is the case we study in this paper.

In particular, we consider a one-dimensional It\^o semimartingale $X$
which is
defined on some probability space $(\Omega,\f,(\f_t)_{t\geq0},\PP
)$ and can always be
represented as
%
%
\begin{eqnarray}
\label{1-1} %
\hspace*{18pt} X_t&=&X_0+
\int_0^tb_s \,ds+\int
_0^t\si_s \,dW_s+
\int_0^t\int_\R
\delta(s,z) 1_{\{|\delta(s,z)|\leq1\}} (\um-\un ) (ds,dz)
\nonumber
\\[-8pt]
\\[-8pt]
&&{}+\int_0^t\int_\R
\delta(s,z) 1_{\{|\delta(s,z)|>1\}} \um (ds,dz),
\nonumber
\end{eqnarray}
where $W$ is a standard Brownian motion and $\um$ a Poisson random
measure on $\R_+\times\R$ with compensator (intensity measure)
$\un(dt,dz)=dt\otimes dz$. This is the Grigelionis representation, and the
specific choice of the Poisson measure $\um$ is in no way a
restriction (see,
e.g., Theorem~2.1.2 in \cite{JP}).
Here, $b$ and $c$ are progressively measurable processes and $\delta$ is
a predictable function on $\Om\times\R_+\times\R$, with appropriate
integrability assumptions.

The process $X$ is observed at regularly spaced times $i\De_n$ for
$i=0,1,\ldots$\,, within a finite time interval $[0,T]$, and without
microstructure noise. Our goal is to estimate, on the basis of these
observations, the so-called integrated volatility, that is,
%
%
\begin{equation}
\label{1-2} C_t=\int_0^tc_s
\,ds \qquad\mbox{where } c_s=\si_s^2,
\end{equation}
for $t=T$ or more generally for all $t\in(0,T]$, with the rate $1/\rdn$,
when $X$ contains jumps of infinite variation.

When jumps are absent, that is, when $\delta\equiv0$ [so the last two
terms in
\eqref{1-1} disappear], the best estimator of $C_t$ is the \emph{realized
volatility}, or \emph{approximate quadratic variation}:
%
%
\begin{equation}
\label{1-3} \wC^n_t=\sum_{i=1}^{[t/\De_n]}
\bigl(\De^n_iX \bigr)^2\qquad\mbox{where
} \De^n_iX=X_{i\De_n}-X_{(i-1)\De_n}.
\end{equation}
Under very weak assumptions on $b$ and $c$ (namely when $\int_0^tb^2_s
\,ds$
and $\int_0^tc_s^2 \,ds$ are finite for all $t$), we have a central limit\vspace*{-3pt}
theorem (CLT) with rate $\frac{1}{\rdn}$: the processes $\frac
{1}{\rdn}(\wC^n_t-C_t)$
converge in the sense of stable convergence in law for processes,
to a limit $Z$ which is defined on an extension of the space, and which
conditionally on $\f$ is a centered Gaussian martingale whose conditional
law is characterized by its (conditional) variance
%
%
\begin{equation}
\label{1-4} V_t := \E \bigl((Z_t)^2\mid\f
\bigr) = 2\int_0^tc_s^2
\,ds,
\end{equation}
or equivalently, we have $Z_t=\sqrt{2}\int_0^tc_s \,dW^{(1)}_s$,
where $W^{(1)}$
is a Brownian motion independent of $\f$. Furthermore, when
$c_s(\omega)=c$ is
a constant, or more generally when $c_t(\omega)=c(t,X_t(\omega))$ for
a smooth
enough function $c$ on $\R_+\times\R$, the estimators $\wC^n_t$ are
\emph{efficient} for any fixed time $t$, because in this case we have
the LAN\vadjust{\goodbreak} or
LAMN property and $V_t$ above is the inverse of the
$\f$-conditional Fisher information, normalized by $\De_n$.
Therefore, in the general case \eqref{1-1} with $\delta\equiv0$ we
qualify the
estimator $\wC^n_t$ as being \emph{efficient}.

When jumps are present, so far there are essentially two types of results,
hinging on a specification of the so-called \emph{degree of jump activity}.
To keep things simple in this \hyperref[sec-Intro]{Introduction}, and
although substantial
extensions can be made, we will suppose that for some $r\in[0,2]$,
%
%
\begin{eqnarray}
\label{1-5} \bigl\llvert \delta(\omega,t,z)\bigr\rrvert ^r\wedge1
\leq J(z)
\nonumber
\\[-8pt]
\\[-8pt]
\eqntext{\mbox{where $J$ is a Lebesgue-integrable function on $\R$}.}
\end{eqnarray}
The smaller $r$ above is, the stronger the assumption is, and it is
(slightly) stronger than assuming $\sum_{s\leq t}|\De X_s|^r<\infty$
for all $t$, where $\De X_s=X_s-X_{s-}$ is the size of the jump at time
$s$. When (\ref{1-5}) holds with $r=0$, the jumps have finite
activity; when
(\ref{1-5}) holds with $r=1$, the jumps are (locally) summable. In the latter
case, we can rewrite \eqref{1-1} (up to modifying $b_t$) as
%
%
\begin{equation}
\label{1-6} X_t=X_0+\int_0^tb_s
\,ds+\int_0^t\si_s
\,dW_s+ \int_0^t\int
_\R\delta(s,z) \um(ds,dz).
\end{equation}
The supremum of all $r$ for which \eqref{1-5} holds is the degree of jump
activity, or Blumenthal--Getoor index. Then we have two cases:
\begin{longlist}
\item[1.]\textit{When $r< 1$}. In this case, we have two major types
of volatility
estimators that enjoy a feasible CLT. The first is the \emph{truncated realized
volatility} (cf. \cite{M1,M2,JP})
%
%
\begin{equation}
\label{1-7} \operatorname{TC}(v_n)^n_t=
\sum_{i=1}^{[t/\De_n]} \bigl(\De^n_iX
\bigr)^2 1_{\{
|\De^n_iX|\leq v_n\}},\qquad v_n\asymp
\De_n^\vpi
\end{equation}
[the last statement means that $\frac{1}A\leq v_n/\De_n^\vpi\leq A$
for some
$A\in(1,\infty)$]. $\operatorname{TC}(v_n)^n_t$ has exactly the same limiting
properties as
$\wC^n$ does in the continuous case provided \eqref{1-5} holds with some
$r\in[0,1)$ and $\vpi\in [\frac{1}{2(2-r)},\frac{1}2 )$.

The second type of jump-robust volatility estimators are the \emph{multipower
variations} (cf. \cite{BS,BSW,JP}),
which we do not explicitly recall here. These estimators also satisfy a
CLT with
rate $\frac{1}{\rdn}$, but with a conditional variance bigger than in
\eqref
{1-4} (so
they are rate-efficient but not variance-efficient).

\item[2.]\textit{When $r\geq1$}. In this case, the above two types
of estimators
are still consistent, but when centered around $C_t$ and appropriately scaled,
they are only bounded in probability with no CLT in general and rate of
convergence that is much slower than $1/\sqrt{\Delta_n}$. For
example, when
$r\geq1$, the sequence $\frac{1}{\De_n^{\vpi(2-r)}}
(\operatorname{TC}(v_n)^n_t-C_t)$ is bounded
in probability (when $r=1$, the multipower variations enjoy a CLT with
a bias
term, see \cite{V}).

On a more general level, we have
the following general result from \cite{JR}: If we have estimators
$\wC^{\prime n}_t$
such that, for some sequence $w_n\to\infty$ of numbers, the variables
$w_n (\wC^{\prime n}_t-C_t)$ are bounded in probability in $n$ and also
when $X$ ranges through all semimartingales of type \eqref{1-1} satisfying
\eqref{1-5} with a fixed function $J$ and also $|b_t|+c_t\leq A$ for some
constant $A$ (so $w_n$ is a kind of ``minimax'' rate), we necessarily
have for
some constant $K$:
%
%
\begin{equation}
\label{1-8} w_n \leq\cases{ %
K/\rdn&\quad$
\mbox{if } 0\leq r\leq1$,
\cr
K \biggl(\displaystyle\frac{\log(1/\De_n)}{\De_n}
\biggr)^{{(2-r)}/{2}}&\quad$\mbox{if } 1<r<2$. 
}
\end{equation}
\end{longlist}

In this paper, we exhibit new estimators for $C_t$ which converge
with rate $\frac{1}{\rdn}$, and which are even variance-efficient in
the sense that
they satisfy the same CLT as $\wC^n_t$ does in the continuous case, when
$r$ defined in (\ref{1-5}) above, that is, the jump activity, is
bigger than $1$.
Of course, given the result in \cite{JR}, discussed in point (2)
above, this
is only possible under some additional assumption, namely that the ``small''
jumps behave like those of a stable process, or of the integral with
respect to a stable-like process, with some index $\be\in(1,2)$
[recall that
in this case \eqref{1-5} holds for all $r>\be$, but not for $r\leq
\be$].
Hence, here we are working in a kind of semiparametric setting, with the
(unknown) parameter $\be$. We should point
out that this ``semiparametric'' setting is still quite general and covers
many jump models used in empirical applications, particularly those in
finance. Similar assumptions about the jumps have been also made when
estimating the Blumenthal--Getoor index of jump activity in \cite{AJ}
and \cite{TT-0} among others.

The estimation method proposed in the current paper is based on estimating
locally the volatility (diffusion coefficient) from the empirical
characteristic function of the increments of the process over blocks of
decreasing length but containing an increasing number of observations, and
then summing the local volatility estimates. The separation of
volatility from
jumps in our method is due to the dominant role of the diffusion
component of
$X$ in (the real part of) the characteristic function of the high frequency
increments of the process for values of the argument that are going to
infinity at the rate $1/\rdn$, or at a slightly slower rate.

When infinite variation jumps are present, the proposed volatility estimators
contain a bias which determines their rate of convergence. The bias scales
differently for different values of the argument of the empirical
characteristic
function, used in forming our nonparametric volatility estimators, and
we use
this property to debias our initial volatility estimators. The debiased
volatility estimators achieve the efficient rate of convergence and
some of
them reach the same (efficient) asymptotic variance as in (\ref{1-4}).

The empirical characteristic function of high-frequency increments has been
previously used in nonparametric estimation of the empirical Laplace
transform of volatility in \cite{TT} as well as in \cite{TT-1} for estimation
of the empirical Laplace transform of the stochastic scale for pure-jump
semimartingales. There are two major differences between these papers and
our study. First, we are interested in estimating the integrated volatility
while the above cited papers consider estimation of the empirical Laplace
transform of the stochastic volatility. Second, and more importantly,
\cite{TT} consider jump-diffusion setting with jumps of finite variation
only and \cite{TT-1} consider pure-jump semimartingales (i.e., processes
with no diffusion). Our main contribution is rate and variance efficient
estimators of integrated volatility in jump-diffusion setting with
jumps of
infinite variation. Finally, the empirical characteristic function in low
frequency setting has been used in \cite{NR,KR} and \cite{R} for
estimating the diffusion coefficient of a L\'{e}vy process, in \cite{CDH}
for nonparametric estimation for a L\'{e}vy process which is a sum of a
drift, a symmetric stable process and a compound Poisson process, as
well as
in \cite{Be} and \cite{BeP} for estimation of L\'{e}vy density and jump
activity in affine models.

The paper is organized as follows. In Section~\ref{sec-Set}, we
present the setting and state our assumptions. In
Section~\ref{sec-FE}, we propose our initial estimators of
integrated volatility and derive a CLT for them when a bias due to
the infinite variation jumps is removed from the estimators. In
Section~\ref{sec-EF}, we propose a way to estimate this bias and
derive a feasible CLT for our debiased estimators.
Section~\ref{sec-MC} contains a Monte Carlo study. Proofs are
given in Section~\ref{sec-PFE}.

\section{The setting}\label{sec-Set}

As mentioned before, the underlying process $X$ is a one-dimensional
It\^o
semimartingale on the space $(\Omega,\f,(\f_t)_{t\geq0},\PP)$, and
observed without noise at the
times $i\De_n$: $i=0,1,\ldots$\,. We restrict the general form \eqref{1-1}
by assuming that the jumps are a mixture of (essentially unspecified) jumps
with finite variation, plus the jumps of a stochastic integral with respect
to a L\'evy process whose small jumps are ``stable-like.''

We have two versions, the simplest one being as follows:
%
%
\begin{equation}
\label{2-1} \qquad X_t=X_0+\int_0^tb_s
\,ds+\int_0^t\si_s
\,dW_s+\int_0^t
\gamma_{s-} \,dY_s+ \int_0^t
\int_\R\delta(s,z) \um(ds,dz)
\end{equation}
with $Y$ a symmetric pure jump L\'evy process with Blumenthal--Getoor index
$\be\in[0,2)$ and the last integral being with finite variation (the precise
assumptions are given below). In this version, the jumps due to $Y$ are
``symmetric'' in the sense that $\int_0^t\gamma_{s-} \,dY_s$ and
$-\int_0^t\gamma_{s-} \,dY_s$ have the same law, as processes. To deal
with the
nonsymmetric case, one could use a process $Y$ which is nonsymmetric.
However, it is more convenient and also more general to use the following
version:
%
%
\begin{eqnarray}
\label{2-2} X_t&=&X_0+\int_0^tb_s
\,ds+\int_0^t\si_s
\,dW_s+ \int_0^t \bigl(
\gamma^{+}_{s-} \,dY^{+}_s+
\gamma^{-}_{s-} \,dY^{-}_s \bigr)
\nonumber
\\[-8pt]
\\[-8pt]
&&{}+
\int_0^t\int_\R
\delta(s,z) \um(ds,dz)\nonumber
\end{eqnarray}
with $Y^+$ and $Y^-$ two independent L\'evy processes with the same
index $\be$ and positive jumps.

We will also require the volatility
$\si_t$ to be an It\^o semimartingale, and it can thus be represented as
%
%
\begin{eqnarray}
\label{2-3} %
\si_t&=&\si_0+
\int_0^tb^\si_s \,ds+
\int_0^tH^\si_s
\,dW_s+\int_0^tH^{\prime \si}_s
\,dW_s'
\nonumber
\\
&&{}+\int_0^t
\int_\R\delta^\si(s,z) 1_{\{|\delta^\si(s,z)|\leq1\}} (\um-
\un) (ds,dz)
\\
&&{}+\int_0^t\int_\R
\delta^\si(s,z) 1_{\{|\delta^\si(s,z)|>1\}} \um(ds,dz).
\nonumber
\end{eqnarray}
Most volatility models used in empirical applications satisfy (\ref{2-3}),
in particular, models in the popular affine class.

As is well known, the jumps of $\si_t$ can, without restriction, be driven
by the same Poisson measure $\um$ as $X$, but we need a second Brownian
motion $W'$: in the case of ``pure leverage,'' we would have $H^{\prime\si}
\equiv0$
and $W'$ is not needed; in the case of ``no leverage,'' we rather have
$H^\si\equiv0$ and in the mixed case we need both $W$ and~$W'$.

Note that \eqref{2-1} is a special case of \eqref{2-2}: indeed, if $Y$
is a pure jump symmetric L\'evy process, it can always be written as
$Y=Y^+-Y^-$ with $Y^+$ and $Y^-$ being independent identically distributed
and with positive jumps, so \eqref{2-2} with $\gamma^+=\gamma$ and
$\gamma^-=-\gamma$
is the same as \eqref{2-1} with $\gamma$. Therefore, we only give the
assumptions
for \eqref{2-2}. The first assumption is a structural assumption
describing the driving terms $W,W',\um,Y^{\pm}$, the second one being
a set of conditions on the coefficients implying in particular the existence
of the various stochastic integrals involved above. Both assumptions involve
a number $r$ in $[0,1)$ (the same in both) and, the smaller $r$ is, the
stronger the two assumptions are.

\renewcommand{\theAssumption}{(A)}
\begin{Assumption}\label{ass(A)}
The processes $W$ and $W'$ are two independent
Brownian motions, independent of $(\um,Y^{+},Y^{-})$;
the measure $\um$ is a Poisson random measure on $\R_+\times\R$
with intensity
$\un(dt,dz)=dt\otimes dz$; the processes $Y^{\pm}$ are two
independent L\'evy\vadjust{\goodbreak}
processes with characteristics $(0,0,F^{\pm})$ and positive jumps
[i.e.,
each $F^{\pm}$ is supported by $(0,\infty)$].
Moreover, there is a number $\be\in[1,2)$ such that
the tail functions $\BF^{\pm}(x)=F^{\pm}((x,\infty))$ satisfy
%
%
\begin{equation}
\label{2-5} x\in(0,1] \quad\Rightarrow\quad \biggl|\BF^{\pm}(x)-
\frac{1}{x^{\be}} \biggr|\leq g(x),
\end{equation}
where $g$ is a decreasing function such that
$\int_0^1x^{r-1} g(x) \,dx<\infty$.
\end{Assumption}

\renewcommand{\theAssumption}{(B)}
\begin{Assumption}\label{ass(B)}
We have a sequence $\tau_n$ of stopping times
increasing to infinity, a sequence $a_n$ of numbers, and a nonnegative
Lebesgue-integrable function $J$ on $\R$, such that
the processes $b,H^\si,\gamma^\pm$ are c\`adl\`ag
adapted, the coefficients $\delta,\delta^\si$ are predictable, the processes
$b^\si,H^{\prime\si}$ are progressively measurable, and
%
%
\begin{eqnarray}
\label{2-6} %
&\hspace*{13pt}t<\tau_n \quad\Rightarrow
\quad\bigl |\delta(t,z)\bigr |^r\wedge1\leq a_nJ(z),\qquad \bigl |
\delta^\si(t,z)\bigr |^2\wedge1\leq a_nJ(z),&
\nonumber
\\
&\hspace*{13pt}t<\tau_n,\qquad V=b,b^\si,H^\si,H^{\prime\si},
\gamma^{+},\gamma^{-} \quad\Rightarrow\quad
|V_t|\leq a_n,&
\nonumber
\\[-12pt]
\\[-4pt]
&\hspace*{13pt}V=b,H^\si,\gamma^{+},\gamma^{-}\hspace*{241pt}&\nonumber
\\
&\hspace*{13pt}\qquad\quad\Rightarrow\quad\bigl |\E(V_{(t+s)\wedge\tau_n}-V_{t\wedge\tau_n}\mid\f_t)\bigr |+ \E
\bigl(\bigl |V_{(t+s)\wedge\tau_n}-V_{t\wedge\tau_n}\bigr |^2\mid\f_t
\bigr)\leq a_ns.&
\nonumber
\end{eqnarray}
\end{Assumption}

Note that we do \emph{not} require the processes $Y^{\pm}$ to be independent
from the measure $\um$, thus allowing any kind of dependence between the
jumps of $X$ and those of $\si$. Intuitively, the number $r$ in
Assumptions~\ref{ass(A)}
and~\ref{ass(B)} control the activity of the finite jump variation
component of
$X$ as well as the degree of deviation from the stable process of
$Y^{\pm}$
which drive the infinite jump variation component of $X$. Our condition in
(\ref{2-5}) is similar to condition AN1 on the L\'{e}vy measure around zero
in \cite{BeP}. Assumptions~\ref{ass(A)} and~\ref{ass(B)} are
satisfied by many parametric
models for the jump component used in applications as illustrated by the
following example.

\begin{Example*}
Suppose the jump component of $X$ is given by a
time-changed L\'evy process with absolute continuous time-change, that is,
$L_{T_t}$ where $L_v$ is a pure-jump L\'evy process with L\'evy measure $F$
satisfying (\ref{2-5}) and time-change $T_t = \int_0^ta_s\,ds$ for
$a_t$ being
strictly positive It\^o semimartingale. A popular parametric example
for $F$
is that of a tempered stable process with corresponding L\'evy density of
the form
\[
\frac{A^{+}e^{-\lambda^{+}x}}{|x|^{1+\beta}}1_{\{x>0\}}+\frac{A^{-}
e^{-\lambda^{-}|x|}}{|x|^{1+\beta}}1_{\{x<0\}},\qquad
A^{\pm}\geq0, \lambda^{\pm}>0, \beta\in(0,2).
\]
In this case, it is not hard to show (using Theorem~2.1.2 of \cite{JP} which
links integrals of random functions with respect
to Poisson measure and random integer-valued measures) that
Assumptions~\ref{ass(A)}
and~\ref{ass(B)} (regarding the
jump part of $X$) hold with $\beta$ in Assumption~\ref{ass(A)}\vadjust{\goodbreak} being
the corresponding
parameter in the above parametric model when $\beta\in[1,2)$ and further
$r=\beta-1+\iota$ for $\iota>0$ arbitrary\vspace*{-1pt} small and $\gamma_t^{+} =
(\frac{A^+a_t}{\beta} )^{1/\beta}$ and $\gamma_t^{-}=-
(\frac{A^-a_t}{\beta} )^{1/\beta}$ [and nonzero $\delta
$ in
(\ref{2-2}) which\vspace*{-1pt} depends on $Y^{\pm}$]. When $\beta\in(0,1)$ in
the above
parametric model, Assumptions~\ref{ass(A)} and~\ref{ass(B)}
hold trivially with $\gamma_t^{\pm} = 0$.
\end{Example*}
%

We end this section with a few comments:
\begin{longlist}
\item[1.] In \eqref{2-5}, there is an implicit standardization of the
processes
$Y^{\pm}$. One could replace it by
%
%
\begin{equation}
\label{0-2B} x\in(0,1] \quad\Rightarrow\quad \biggl|\BF^{\pm}(x)-
\frac{a_{\pm}}{x^{\be}} \biggr|\leq g(x)
\end{equation}
for positive constants $a_{\pm}$. However, in this case the
processes\vspace*{-1pt} $Y^{\prime\pm}=Y^{\pm}/a_{\pm}^{1/\be}$ satisfy \eqref{2-5} as
stated, and \eqref{2-2} holds with $Y^{\prime\pm}$ and $\gamma^{\prime\pm}=
a_{\pm}^{1/\be} \gamma^{\pm}$ as well. It is more convenient in the sequel,
and not a restriction, to use the standardized form \eqref{2-5}.

\item[2.] By Assumption~\ref{ass(B)} and the fact that $r<1$, the
last integral in
\eqref{2-2} defines a process with finite variation which is the sum of
its jumps. On the other hand, $\int_0^t(\gamma^{+}_{s-} \,dY^{+}_s
+\gamma^{-}_{s-} \,dY^{-}_s)$
has a Blumenthal--Getoor (BG) index $\be\geq1$ and is of infinite variation
(even when $\be=1$, and unless $\gamma^{+}$ and $\gamma^{-}$ identically
vanish, of course), although still a (compensated) ``pure jump'' process.

\item[3.] Concerning the regularity assumptions in \ref{ass(B)}, the
last part of
\eqref{2-6} could be somewhat weakened (e.g., we could drop it
in the
case of $V=H^\si$), but at the expense of a nontrivial
complication of the proofs. Since these are satisfied in virtually all
models used in practice, we decided to impose these assumptions here.
Note also that this last part of \eqref{2-6} is satisfied as soon as the
processes $b,H^\si,\gamma^+,\gamma^-$ are themselves It\^o semimartingales
with locally bounded characteristics.\vspace*{-1pt}
\end{longlist}

\section{First estimators of $C_t$}\label{sec-FE}

In this section, we construct our initial estimators of $C_t$. These
estimators are not efficient in general, but they will be used
to construct efficient estimators later on.

We use the real part of the ``local'' (in time) empirical characteristic
functions of increments, taken at point $u_n/\rdn$ for some sequence
$u_n>0$ going to $0$ slowly enough. Here, ``local'' means
that the empirical characteristic function is constructed on windows of time
length $v_n$ or $2v_n$, where
$v_n=k_n\De_n$ and $k_n\geq1$ is a suitable sequence of integers
going to
infinity, to be specified later. We will in fact use two different versions:\vspace*{-1pt}
%
%
\begin{eqnarray}
\label{F-1} %
&&\mbox{Symmetrized version:}\nonumber
\\[-1pt]
&&\qquad
L(u)^n_j = \frac{1}{k_n} \sum
_{l=0}^{k_n-1} \cos \bigl(u \bigl(\De^n_{2jk_n+1+2l}X-
\De^n_{2jk_n+2+2l}X \bigr)/\rdn \bigr),
\nonumber
\\[-8.5pt]
\\[-8.5pt]
&&\mbox{Nonsymmetrized version:}\nonumber
\\[-1pt]
&&\qquad L'(u)^n_j =
\frac{1}{k_n} \sum_{l=0}^{k_n-1}\cos
\bigl(u\De^n_{1+jk_n+l}X/\rdn \bigr)
\nonumber
\end{eqnarray}
for $j\geq1$ some integer, $u>0$ some real and recall
$\De^n_iX = X_{i\De_n}-X_{(i-1)\De_n}$.
$L(u)^n_j$ and $L'(u)^n_j$ are not bigger than $1$, and the variables\vspace*{-1pt}
%
%
\begin{eqnarray}
\label{F-2} %
\widehat{c}(u)^n_j&=&-
\frac{1}{u^2} \log \biggl(L(u)^n_j\vee
\frac{1}{\sqrt{k_n}} \biggr),
\nonumber
\\[-8.5pt]
\\[-8.5pt]
\widehat{c}'(u)^n_j&=&-
\frac{2}{u^2} \log \biggl(L'(u)^n_j
\vee\frac{1}{\sqrt{k_n}} \biggr),
\nonumber
\end{eqnarray}
satisfy $0\leq\widehat{c}(u)^n_j\leq\frac{\log k_n}{2u^2}$ and
$0\leq\widehat{c}'(u)^n_j\leq\frac{\log k_n}{u^2}$, and serve as
local\vspace*{1pt} estimators of the volatility\vadjust{\goodbreak} (of the average of $c_t$ over the
interval $(2jv_n,2(j+1)v_n]$ or $(jv_n,(j+1)v_n]$, to be more precise). The
associated estimators for integrated volatility are thus (recall
$v_n=k_n\De_n$):
%
%
\begin{eqnarray}
\label{F-3} %
\wC(u)^n_t&=&2v_n
\sum_{j=0}^{[t/2v_n]-1} \biggl(\widehat
{c}(u)^n_j-\frac{1}{u^2k_n} \bigl(\sinh
\bigl(u^2\widehat{c}(u)^n_j \bigr)
\bigr)^2 \biggr),
\nonumber
\\[-8pt]
\\[-8pt]
\wC'(u)^n_t&=&v_n\sum
_{j=0}^{[t/v_n]-1} \biggl(\widehat {c}'(u)^n_j-
\frac{2}{u^2k_n} \bigl(\sinh \bigl(u^2\widehat{c}'(u)^n_j/2
\bigr) \bigr)^2 \biggr),
\nonumber
\end{eqnarray}
where recall $\sinh(x) = \frac{e^x-e^{-x}}{2}$. On an intuitive level,
$\wC(u)^n_t$ and $\wC'(u)^n_t$ separate volatility (of the
diffusive part of $X$) from jumps in $X$ by utilizing the fact that the
diffusive component of $X$ dominates the behavior of the real part of the
empirical characteristic function at high-frequencies for values of the
argument
that are ``sufficiently'' away from zero. Indeed, in the simple case when
$X_t = X_0 + bt+\sigma W_t+\gamma Y_t$ for $Y_t$ a symmetric $\beta$-stable
process with unit scale, we have $\log\Re(\mathbb{E}(e^{iu\Delta_i^nX
/\sqrt{\Delta_n}})) = \log(\cos(ub\Delta_n^{1/2}))-\frac
{u^2\sigma^2}{2}-
|\gamma|^{\beta}u^{\beta}\Delta_n^{1-\beta/2}$.

The terms $\frac{1}{u^2k_n}  (\sinh(u^2\widehat{c}(u)^n_j)
)^2$ and
$\frac{2}{u^2k_n}  (\sinh(u^2\widehat{c}'(u)^n_j/2) )^2$
remove biases of
higher asymptotic order in $\widehat{c}(u)^n_j$ and $\widehat
{c}'(u)^n_j$, respectively,
which arise due to the nonlinear transformation of $L(u)^n_j$ and $L'(u)^n_j$
in forming $\widehat{c}(u)^n_j$ and $\widehat{c}'(u)^n_j$.\vspace*{-2pt}

We note that for any fixed $n$, $\lim_{u\downarrow0}\wC'(u)^n_t =
\sum_{i=1}^{[t/\De_n]}(\Delta_i^nX)^2$ is the realized volatility
(which in
presence of jumps does not estimate the integrated volatility). The robustness
of our estimator $\wC'(u)^n_t$ with respect to jumps in $X$ will
result from
using $u=u_n$ that is ``sufficiently'' far from zero, and the
variance-efficiency of the corrected second-step estimators will come from
the fact that $u_n\to0$ (we make this formal in the theorems below).

For stating the asymptotic behavior of the estimators in (\ref{F-3}),
we need
some additional notation. First, for $\be\in(0,2)$ we set
%
%
\begin{eqnarray}
\label{F-4} \be>1 \quad&\mapsto&\quad\chi'(\be)=\int
_0^\infty\frac{1-\cos
y}{y^\be}\, dy,
\nonumber
\\[-8pt]
\\[-8pt]
 \be>0
\quad&\mapsto&\quad\chi(\be)=-\be\chi'(\be+1)= \int
_0^\infty\frac{\sin y}{y^\be} \,dy\nonumber
\end{eqnarray}
(the last integral is convergent for all $\be>0$, but absolutely convergent
when $\be>1$ only). Next, with the notation $\{x\}^\be=|x|^\be\sign
(x)$ for
any $x\in\R$, we associate with the processes $\gamma^{\pm}$
the following [when $\chi'(\be)$ appears below we implicitly suppose
$\be>1$]:
%
%
\begin{eqnarray}
\label{F-5} %
a_t&=&\chi(\be) \bigl(\bigl |
\gamma^{+}_t\bigr |^{\be}+\bigl |\gamma^{-}_t\bigr |^{\be}
\bigr),\qquad a_t'=\chi'(\be) \bigl(
\bigl\{\gamma^{+}_t \bigr\}^{\be}+ \bigl\{
\gamma^{-}_t \bigr\}^{\be
} \bigr),
\nonumber
\\
A(u)_t^n&=&2u^{\be-2}\De_n^{1-\be/2}
\int_0^ta_s \,ds,
\\
A'(u)_t^n&=& \frac{2}{u^2}\int
_0^t \bigl(\De_n^{1-\be/2}u^{\be}a_s-
\log \bigl(\cos \bigl(\De_n^{1-\be/2}u^{\be}a_s'
\bigr) \bigr) \bigr) \, ds.
\nonumber
\end{eqnarray}

Under appropriate assumptions on the sequence $u_n$, we will see that
$\wC(u_n)_T$ and $\wC'(u_n)_T$ converge to $C_T$, and there is an associated
central limit theorem with the convergence rate $1/\rdn$. However, in the
CLT there is typically a nonnegligible bias due to the infinite variation
jumps in $X$, and to account for
this bias we consider the following normalized error processes:
%
%
\begin{eqnarray}
\label{F-6} %
Z(u)^n_t&=&
\frac{1}{\rdn} \bigl(\wC(u)^n_t-C_t-A(u)^n_t
\bigr),
\nonumber
\\[-8pt]
\\[-8pt]
Z'(u)^n_t&=&\frac{1}{\rdn} \bigl(
\wC'(u)^n_t-C_t-A'(u)^n_t
\bigr).
\nonumber
\end{eqnarray}

$A(u)^n_t$ and $A'(u)^n_t$ are easiest to understand in the L\'evy
case, that is,
when $\gamma_t^{\pm}$ are constants. In this case, $A'(u)^n_1$ is
$-\frac{2}{u^2}$ times the logarithm of the real part of the characteristic
function of $\De^n_iL/\rdn$, where $L=\gamma^+L^++\gamma^-L^-$ and\vspace*{-1pt}
$L^+$ and $L^-$ are two independent one-sided stable processes with L\'evy
density $\frac{\beta}{x^{\beta+1}} 1_{\{x>0\}}$, and $A(u)^n_1=A'(u)^n_1$
when $\gamma^-=-\gamma^+$. In this case of constant $\gamma^\pm_t$,
taking the difference $\De^n_{i+1}X-\De^n_iX$ makes the contribution of
the stochastic integrals w.r.t. $Y^\pm$ globally symmetric: the characteristic
function of $\De^n_{i+1}L-\De^n_iL$ above becomes real, and this is
why we put $A(u)^n_t$ instead of $A'(u)^n_t$ in the first case of
\eqref{F-6}. Now, $A(u)^n_t$ has a much simpler form than $A'(u)^n_t$,
regarding its dependence upon $u$, which makes its estimation from the data,
as conducted in the next section, rather easy. On the other hand, differencing
increments results in a loss of information, since in the definition of
$\wC(u)^n_t$ we have twice less summands than in the definition of
$\wC'(u)^n_t$. [Note that the form (\ref{2-1}) for $X$ corresponds to having
$\gamma^-=-\gamma^+$, hence in this case $A'(u)^n_t=A(u)^n_t$.]

In order to give a simple version of the limits below, we consider
an extension $(\widetilde{\Omega},\widetilde{\f},(\widetilde{\f
}_t)_{t\geq0},\widetilde{\PP})$ of the original space $(\Omega,\f
,(\f_t)_{t\geq0},\PP)$, which supports two
independent Brownian motions $W^{(1)}$ and $W^{(2)}$, independent of the
$\si$-field $\f$, and on this extension we introduce the two processes
%
%
\begin{equation}
\label{F-12bA} Z_t=\sqrt{2} \int_0^tc_s
\,dW^{(1)}_s,\qquad \BZ_t=\frac{1}{\sqrt{6}}
\int_0^tc_s^2
\,dW^{(2)}_s.
\end{equation}
An equivalent characterization of the pair $(Z,\BZ)$ is as follows:
they are
defined on
an extension $(\widetilde{\Omega},\widetilde{\f},(\widetilde{\f
}_t)_{t\geq0},\widetilde{\PP})$ of $(\Omega,\f,(\f_t)_{t\geq
0},\PP)$ and, conditionally on $\f$, they are
centered continuous Gaussian martingales
characterized by their (conditional) variances--covariances, as given by
%
%
\begin{eqnarray}
\label{F-13b} \WE \bigl((Z_t)^2\mid\f \bigr)&=&2\int
_0^tc_s^2 \,ds,\qquad
\WE \bigl(\BZ_t^2\mid\f \bigr)=\frac{1}6
\int_0^tc_s^4 \,ds,
\nonumber
\\[-8pt]
\\[-8pt]
\WE(Z_t \BZ_t\mid\f)&=&0.\nonumber
\end{eqnarray}
In view of the
debiasing procedure later on, we need a multidimensional version of the CLT,
namely the convergence for all $\theta u_n$, where $\theta$ runs
through a
finite subset $\Te$ of $(0,\infty)$. We are now ready to state the
main results of this section.

%
\begin{theo}\label{TF-1} Assume \ref{ass(A)} and~\ref{ass(B)}
with $r<1$, and choose
$k_n$ and
$u_n$ in such a way that
%
%
\begin{eqnarray}
\label{F-8} k_n \rdn&\to&0,\qquad k_n
\De_n^{1/2-\ep}\to\infty\qquad\forall\ep>0,
\nonumber
\\[-8pt]
\\[-8pt]
 u_n
&\to&0, \qquad\sup_n\frac{k_n\rdn}{u_n^4}<\infty.\nonumber
\end{eqnarray}\vspace*{-\baselineskip}
\begin{longlist}
\item[(a)]We have the (functional) stable convergence in law
%
%
\begin{eqnarray}
\label{F-12}
&&\biggl(Z(u_n)^n, \biggl(\frac{1}{u_n^2}
\bigl(Z(\theta u_n)^n-Z(u_n)^n\bigr)
\biggr)_{\theta\in\Te} \biggr)
\nonumber
\\[-8pt]
\\[-8pt]
&&\qquad\tolls \bigl(\sqrt{2} Z, \bigl(2\sqrt{2} \bigl(
\theta^2-1 \bigr)\BZ \bigr)_{\theta\in\Te} \bigr).\nonumber
\end{eqnarray}
\item[(b)]
If further $\be>1$, we also have
%
%
\begin{equation}
\label{F-12b} \biggl(Z'(u_n)^n, \biggl(
\frac{1}{u_n^2} \bigl (Z'(\theta u_n)^n-Z'(u_n)^n\bigr)
\biggr)_{\theta\in\Te} \biggr) \tolls \bigl(Z, \bigl( \bigl(\theta^2-1
\bigr)\BZ \bigr)_{\theta\in\Te} \bigr).\hspace*{-35pt}
\end{equation}
\end{longlist}
\end{theo}

This exhibits\vspace*{-1pt} a kind of degeneracy. Indeed, \eqref{F-12} and
\eqref{F-12b} imply the following convergence ($\toucp$ means
convergence in
probability, uniformly on each compact time interval):
%
%
\begin{equation}
\label{F-10} Z(\theta u_n)^n-Z(u_n)
\toucp0,\qquad Z'(\theta u_n)^n-Z'(u_n)
\toucp0.
\end{equation}

%
\begin{rema}\label{RR-1}
A possible choice for $k_n$ and $u_n$ is
$k_n\asymp1/\rdn(\log(1/\De_n))^x$ and $u_n\asymp1/(\log(1/\De_n))^{x'}$,
which satisfies \eqref{F-8} as soon as the reals $x,x'$ are such that
$0<x'\leq\frac{x}4$.
\end{rema}

\section{Efficient estimators of $C_t$}\label{sec-EF}

In general, the bias terms $A(u)^n_t$ or $A'(u)^n_t$ in \eqref{F-6} determine
the second-order behavior of the estimators $\wC(u)^n_t$ and
$\wC'(u)^n_t$, thus preventing rate efficiency. In one important
case, though, Theorem~\ref{TF-1} implies that $\wC'(u)^n_t$ will be
both rate
and variance efficient and $\wC(u)^n_t$ will be rate efficient but with
asymptotic variance somewhat larger. This is the
case when the jumps in $X$ are of finite variation, that is, when
$\gamma^{+}$ and
$\gamma^{-}$ are identically $0$. Then \eqref{F-6} reduces to
\[
Z(u)^n_t=\rdn \bigl(\wC(u)^n_t-C_t
\bigr),\qquad Z'(u)^n_t=\rdn \bigl(
\wC'(u)^n_t-C_t \bigr),
\]
and Theorem~\ref{TF-1} implies:

%
\begin{theo}\label{TEF-1} Assume \ref{ass(A)} and~\ref{ass(B)}
with $\gamma^\pm\equiv
0$ and $r<1$,
and choose
$k_n$ and $u_n$ satisfying \eqref{F-8}. Then the processes $Z(u_n)^n$ and
$Z'(u_n)^n$ converge stably in law to $\sqrt{2} Z$ and $Z$, respectively.
\end{theo}

This means, in particular, that the estimators $\wC'(u_n)_t$ are asymptotically
equivalent to the truncated realized volatility $\operatorname
{TC}(v_n)_t$ of \eqref{1-7}
with\vspace*{-1pt} $v_n\asymp\De_n^{\vpi}$ and $\vpi\in (\frac{1}{2(2-r)}\frac{1}2 )$,
and hence are rate and variance efficient. Thus, we provide an
alternative to the truncated realized volatility which is important in
applications due to the presence of tuning parameters in the
construction of both jump-robust volatility estimators (ours and the
truncated realized volatility).

%
\begin{rema}\label{RE-1}
Whereas the above is a special case of Theorem~\ref{TF-1}, it is possible [although far from trivial when one allows the
process $\si$ to jump, as in \eqref{2-3}] to show that when again
$\gamma^\pm\equiv0$ and when $r=1$, and if we fix $u>0$, then the sequence
$Z'(u)^n$ stably converges in law to a process $Z(u)$ which has the same
description as $Z$ above, except that the conditional variance is now
%
%
\begin{equation}
\label{2-9} \WE \bigl(Z(u)_t^2\mid\f \bigr)= 8\int
_0^t \biggl(\frac{\sinh(u^2c_s/2)}{u^2}
\biggr)^2 \,ds
\end{equation}
[when $u_n\to0$ we do not know the behavior of $Z'(u_n)^n$]. Hence, the
estimators $\wC'(u)^n_t$ are still rate efficient, but no longer
variance efficient.\vadjust{\goodbreak} However, the right-hand side of \eqref{2-9} goes to
$2\int_0^tc_s^2 \,ds$ as $u\to0$: so, upon choosing $u$ small enough,
one can
approach variance efficiency as close as one wants to.

Note that, even without variance efficiency, the rate efficiency above plus
the fact that the limit is conditionally unbiased seems to be a new
result when $r=1$.
\end{rema}

When\vspace*{-1pt} the term $\int_0^t(\gamma^{+}_{s-} \,dY^{+}_s+\gamma^{-}_{s-}\,
dY^{-}_s)$
in \eqref{2-2} is present, the
estimators $\wC(u_n)^n_t$ and $\wC'(u_n)^n_t$ converge to $C_t$ at a rate
arbitrarily close to $1/\De_n^{(2-\be)/2}$, which up to a logarithmic term
is in accordance with the minimax rate given in \eqref{1-8} (see \cite{JR}).
However, this does
not give us a feasible limit theorem. In this situation,
one can find a way of eliminating the bias term and come up with estimators
with rate $1/\rdn$ and which are even variance
efficient [of course this is possible under Assumptions~\ref{ass(A)}
and~\ref{ass(B)}
only].

To do this, we fix the time horizon $T>0$, and we set
%
%
\begin{eqnarray}
\label{EF-6} %
\wC(u,\ze)^n_T&=&
\wC(u)^n_T-\frac{(\wC(\ze u)^n_T-\wC(u)^n_T)^2} {
\wC(\ze^2 u)^n_T-2\wC(\ze u)^n_T+\wC(u)^n_T},
\nonumber
\\[-8pt]
\\[-8pt]
\wC'(u,\ze)^n_T&=&\wC'(u)^n_T-
\frac{(\wC'(\ze u)^n_T-\wC'(u)^n_T)^2} {
\wC'(\ze^2 u)^n_T-2\wC'(\ze u)^n_T+\wC'(u)^n_T}.
\nonumber
\end{eqnarray}

The new estimators above are biased-corrected analogues of $\wC
(u)^n_T$ and
$\wC'(u)^n_T$. Our estimation of the bias is very intuitive. It
utilizes the
fact that the only difference (asymptotically) in $\wC(u)^n_T$ and
$\wC'(u)^n_T$ for different values of $u$ stems from the presence of
$A(u)_t^n$ and $A'(u)_t^n$. This suggests an easy way to estimate these
biases from the differences of $\wC(u)^n_T$ and $\wC'(u)^n_T$ over different
values of $u$. The next theorem derives the asymptotic behavior of
$\wC(u,\ze)^n_T$ and $\wC'(u,\ze)^n_T$.

%
\begin{theo}\label{TEF-2} Assume \ref{ass(A)} and~\ref{ass(B)} with
$r<1$ and $C_T>0$ a.s.
Choose $k_n$
and $u_n$ satisfying \eqref{F-8} and any $\ze>1$.
\begin{longlist}
\item[(a)]The variables $\frac{1}{\rdn}(\wC(u_n,\ze)^n_T-C_T)$
converge stably in law\vspace*{-2pt}
to the variable $\sqrt{2} Z_T$, which conditionally on $\f$ is
centered Gaussian with (conditional) variance $4\int_0^Tc_s^2 \,ds$.
\item[(b)] Assume\vspace*{-1pt} further that either $1<\be<\frac{3}2$, or that $\be
\geq
\frac{3}2$
and $\gamma^{+}=-\gamma^{-}$ identically. The variables
$\frac{1}{\rdn}(\wC'(u_n,\ze)^n_T-C_T)$ converge stably in law
to the variable $Z_T$, which conditionally on $\f$ is
centered Gaussian with\vspace*{1pt} (conditional) variance $2\int_0^Tc_s^2 \,ds$.
\end{longlist}
In particular, this applies when \eqref{2-2} reduces to \eqref{2-1},
under the only condition $1<\be<2$.
\end{theo}

The estimator $\wC(u_n,\ze)^n_T$ applies in all cases of
Assumptions~\ref{ass(A)}
and~\ref{ass(B)} and is rate efficient but not variance efficient.
$\wC'(u_n,\ze
)^n_T$ is
both variance and rate efficient and no prior knowledge of $\beta$ is needed
(except that $\beta=1$ is excluded) whenever $\gamma^+=-\gamma^-$
which is
the case in many models. When $\gamma^+\neq-\gamma^-$, then we can use
$\wC'(u_n,\ze)^n_T$ only when $\beta<3/2$.

Alternatively, we could iterate the debiasing procedure and achieve
rate and
variance efficiency even in the asymmetric case $\gamma^+\neq-\gamma^-$.
Such an iteration also permits to replace the fourth term on the
right-hand side
of \eqref{2-2} by a sum of $M$ terms $\int_0^t(\gamma^{m+}_{s-} \,dY^{m+}_s
+\gamma^{m-}_{s-} \,dY^{m-}_s)$, with $Y^{m\pm}$ having Blumenthal--Getoor
indices $\be_m$ with $1\leq\be_M<\cdots<\be_1<2$, and under appropriate
conditions. We leave such extensions for future work.

%
\begin{rema}\label{RR10} When $\PP(C_T>0)<1$, the result as stated
may fail. However, a classical argument shows that it still holds
\emph{in restriction to the set} \mbox{$\{C_T>0\}$}.
\end{rema}

\section{Monte Carlo study}\label{sec-MC}

We test the performance of our new method of estimating integrated volatility
and compare it with that of the truncated realized volatility on
simulated data from
the following stochastic volatility model:
%
%
\begin{eqnarray}
\label{eq:mc_sv} X_t = X_0+\int_0^t
\sqrt{c_s} \,dW_s+\eta Y_t,
\nonumber
\\[-8pt]
\\[-8pt]
\eqntext{c_t =c_0+ \displaystyle\int_0^t0.03(1.0-c_s)
\,ds+0.15\int_0^t\sqrt{c_s}
\,dW_s', \eta\geq0,}
\end{eqnarray}
where $W_t$ and $W_t'$ are two independent Brownian motions and $Y_t$
is a
symmetric $\beta$-stable process independent from $W_t$ and $W_t'$. The
volatility $c_t$
is a square-root diffusion process, which is widely used to model stochastic
volatility in financial applications. The parameters of the volatility
specification are set so that the mean and persistence of volatility is
similar to that in actual financial data. In particular, its mean\vadjust{\goodbreak} is $1$
in the stationary case. Since the key advantage of our estimation
procedure is its ability to recover integrated volatility in presence of
infinite variation jumps, in the Monte Carlo we experiment with values of
the stability parameter of $Y_t$ of $\beta=1.25$, $\beta= 1.50$ and $1.75$.
We further vary the constant $\eta$ (in the interval $[0,2]$) which controls
the relative contribution of $Y_t$ in the total variation of~$X_t$.

In the Monte Carlo, we fix the time span to $1$ day (our unit of time
is a
day) and we consider $1/\Delta_n=2400$ and $1/\Delta_n=4800$, which
corresponds to sampling at $10$ and $5$ seconds, respectively, in a
$6.5$-hour trading day. We set $k_n=240$ for $1/\Delta_n=2400$ and we
increase it to $k_n=320$ when $1/\Delta_n=4800$, which correspond to $10$
and $15$, respectively, blocks per unit of time. Experiments with more blocks
per day led to very similar results.

We test in the Monte Carlo the performance of the bias-corrected estimator
$\widehat{C}'(u,\zeta)^n$ defined in (\ref{EF-6}), whose
implementation we now
discuss. The choice of the tuning parameter $u_n=u$ in
$\widehat{C}'(u,\zeta)^n$ plays a nontrivial role. From the asymptotic
variance in (\ref{2-9}), it is clear that a big or small value of $u$
is always
with respect to the level of volatility $c_s$. For this reason, for each
time interval $[t,t+1]$ we set $u$ for that day to
$u_t^n = \frac{1}{(\log(1/\Delta_n))^{1/30}}\frac{1}{\sqrt{\mathit
{BV}_{[t-1,t]}}}$,
where
%
%
\begin{equation}
\mathit{BV}_{[t-1,t]} = \frac{\pi}{2}\sum
_{i=[(t-1)/\Delta_n]+2} ^{[t/\Delta_n]}\bigl |\Delta_{i-1}^nX\bigr |\bigl |
\Delta_{i}^nX\bigr |
\end{equation}
is the bipower variation on the unit interval $[t-1,t)$ which is a consistent
estimator of $\int_{t-1}^{t}c_s\,ds$ that does not require any choice of
tuning parameters. Our time-varying $u_t^n$ is analogous to the selection
of a time-varying threshold for the truncated realized volatility that
is typically
done (and we implement as well here). The scale\vspace*{-1pt} factor
$\frac{1}{(\log(1/\Delta_n))^{1/30}}$ is chosen so that $u_t^n$ converges
to zero very slowly as $\Delta_n\rightarrow0$.

The bias correction term in $\widehat{C}'(u,\zeta)^n_T$ can be split into
the product of two terms, as $(\wC'(\ze u)^n_T-\wC'(u)^n_T)\times
\frac{(\wC'(\ze u)^n_T-\wC'(u)^n_T)} {
\wC'(\ze^2 u)^n_T-2\wC'(\ze u)^n_T+\wC'(u)^n_T}$. The\vspace*{-2pt} first term
is an
estimator for $A'(u)_T^n$, which is time-varying and the second is an
estimator of $\frac{1}{\zeta^{\beta-2}-1}$ which depends only on the
parameter $\beta$. To reduce the noise in our estimate of the bias,
therefore, we use a horizon of $132$ days (6 months) to estimate the second
term, similar to earlier studies on estimation of the Blumenthal--Getoor
index (\cite{AJ} and \cite{TT-0}), and daily data to estimate the
first term
(as the limit of this term is time-varying). Also for the calculation
of the
second term, we use a smaller value of $u$ as this allows to capture
the slope
of $\widehat{C}'(u,\zeta)^n_T$ better. Overall, for a period of
$T=132$ days,
our daily estimator is
%
%
\begin{eqnarray}
\label{eq:mc_chat} \widehat{C}' \bigl(u_t^n,
\zeta \bigr)^n_{[t,t+1]} &=& \widehat{C}'
\bigl(u_t^n \bigr)^n_{[t,t+1]}
-S_T^n \bigl( \bigl(\wC' \bigl(\ze
u_t^n \bigr)^n_{[t,t+1]}-
\wC' \bigl(u_t^n \bigr)^n_{[t,t+1]}
\bigr) \wedge0 \bigr),
\nonumber
\\
\eqntext{t=1,\ldots,T-1,}
\\
S_T^n &=& \sum_{t=1}^T\bigl(\wC'\bigl(0.3\ze u_t^n\bigr)^n_{[t,t+1]}
-\wC'\bigl(0.3u_t^n\bigr)^n_{[t,t+1]}\bigr)
\\
&&{}\bigg/\sum_{t=1}^T\bigl(\wC'\bigl(0.3\ze^2 u_t^n\bigr)^n_{[t,t+1]}\nonumber
\\
&&\phantom{\bigg/\sum_{t=1}^T(}{}-2\wC'\bigl(0.3\ze u_t^n\bigr)^n_{[t,t+1]}
+\wC'\bigl(0.3u_t^n\bigr)^n_{[t,t+1]}\bigr)
\wedge0.
\nonumber
\end{eqnarray}
The restrictions on the sign above are finite sample restrictions with no
asymptotic effect. In the calculation of the bias correction term, we set
$\ze= 1.5$. Finally, if $\widehat{C}'(u_t^n,\zeta)^n_{[t-1,t]}$ is negative
we repeat the calculation in (\ref{eq:mc_chat}) with $2u_t^n/3$ (this again
has no asymptotic effect).

%
\begin{figure}

\includegraphics{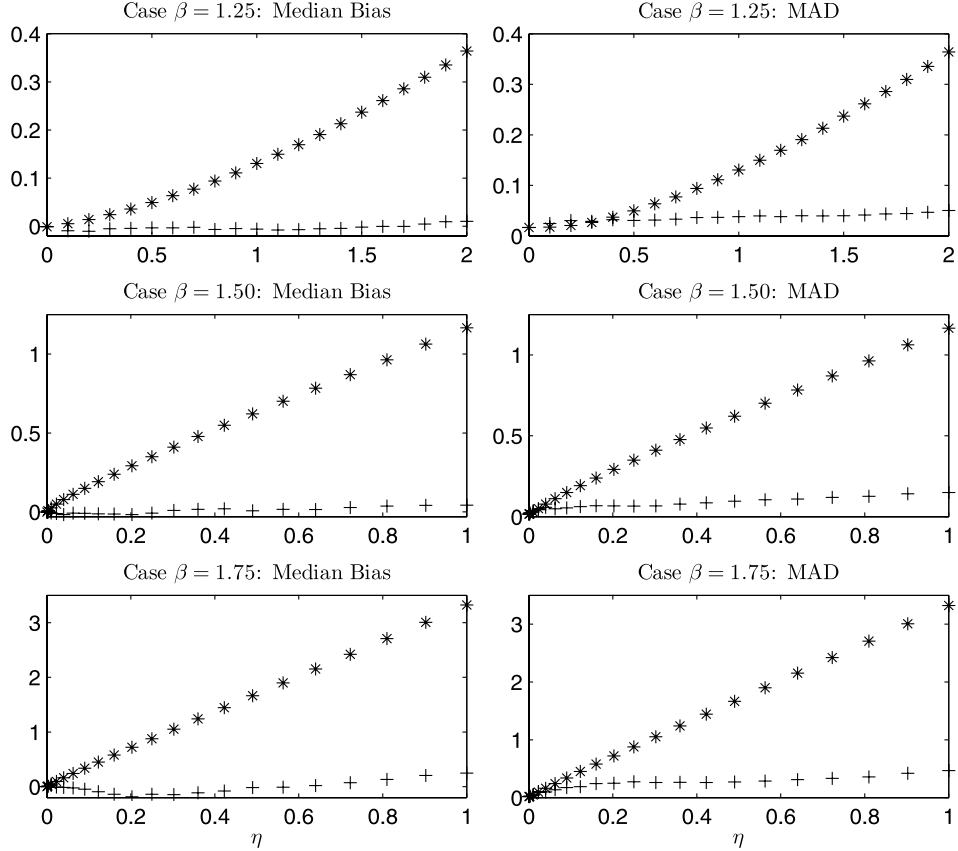}

\caption{Median bias and Median Absolute Deviation (MAD) around the true
value, $\int_t^{t+1}c_s\,ds$, for sampling frequency $1/\Delta_n =
2400$. $+$
corresponds to $\widehat{C}'(u,\zeta)^n$ and $*$ to $\operatorname
{TC}(v_n)^n$.}
\label{fig:mc_bias_mad_sqrt_2400}
\end{figure}

For the truncation realized volatility estimator $\operatorname
{TC}(v_n)^n$, which we
compare below to
our estimator, we set $v_n = 4\sqrt{\mathit{BV}_{[t-1,t]}}\Delta
_n^{0.49}$, as
typically done in existing work.

%
\begin{figure}

\includegraphics{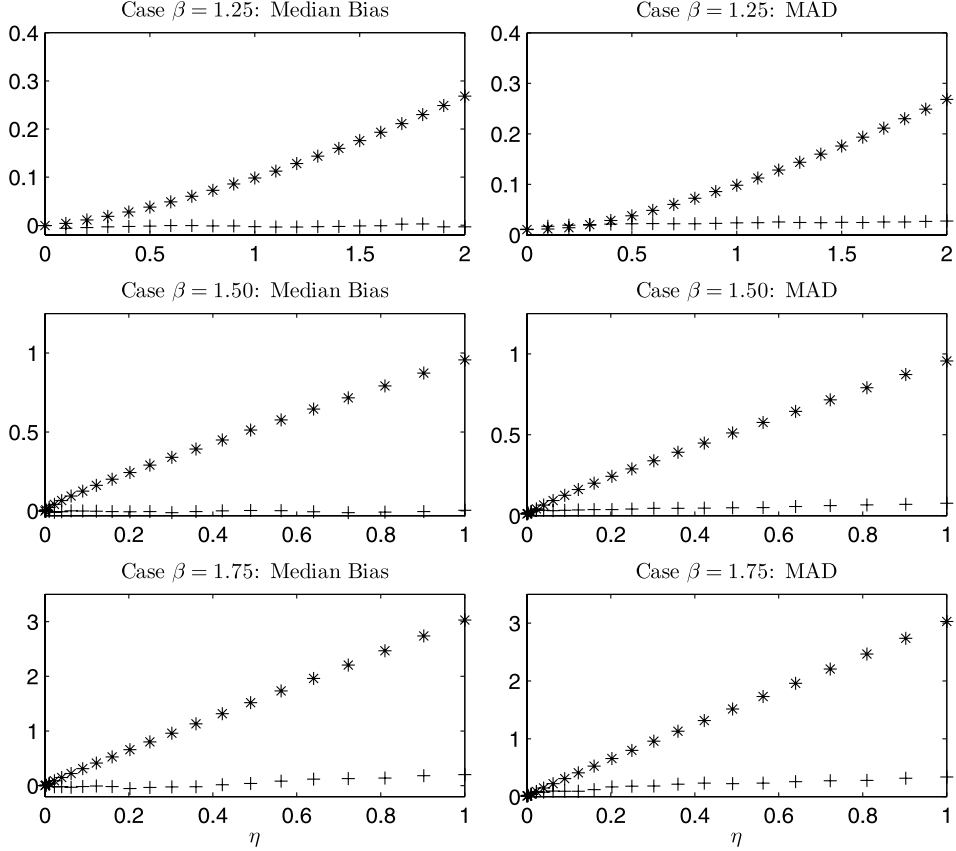}

\caption{Median bias and Median Absolute Deviation (MAD) around the true
value, $\int_t^{t+1}c_s\,ds$, for sampling frequency $1/\Delta_n = 4800$.
$+$ corresponds to $\widehat{C}'(u,\zeta)^n$ and $*$ to
$\operatorname{TC}(v_n)^n$.}
\label{fig:mc_bias_mad_sqrt_4800}
\end{figure}

The results from the Monte Carlo are summarized in
Figures~\ref{fig:mc_bias_mad_sqrt_2400} and~\ref{fig:mc_bias_mad_sqrt_4800}.
Not surprisingly, the activity of the jump component (controlled by
$\beta$)
and its relative share in total return variation (controlled by $\eta
$) have
clear impact on the ability to separate integrated variance from the jumps
in $X$. Our volatility estimator $\widehat{C}'(u_t^n,\zeta)^n_{[t-1,t]}$
performs significantly better than the truncated variance in presence of
infinite variation jumps (recall that both estimators are consistent
regardless of the activity of the jumps). The superior performance of
$\widehat{C}'(u_t^n,\zeta)^n_{[t-1,t]}$ is largely due to the removal of
the bias in the volatility estimation that is due to the infinite variation
jumps. As a result $\widehat{C}'(u_t^n,\zeta)^n_{[t-1,t]}$, unlike
$\operatorname{TC}(v_n)^n$, is essentially unbiased in all considered cases.
Increasing the
sampling frequency improves the performance of both estimators in all cases.
We note, however, that the reduction of bias and MAD for $\operatorname
{TC}(v_n)^n$
for the
higher jump activity case ($\beta=1.75$) is significantly slower and this
is unlike our estimator. This is consistent with the slow rate of convergence
of $\operatorname{TC}(v_n)^n$ in the case of infinite variation jumps
discussed in the
\hyperref[sec-Intro]{Introduction}. Overall, we conclude that our estimator provides a nontrivial
improvement over existing methods for the nonparametric estimation of
integrated volatility in presence of infinite variation jumps.
\eject

\section{Proofs}\label{sec-PFE}

\subsection{Preliminaries}

By a standard localization procedure, we may and will assume that in
\ref{ass(B)}
we have $\tau_1\equiv\infty$ and $J$ is bounded, and also that $X$
and $\si$
are themselves bounded, as well as the jumps of $Y^{\pm}$. Up to
modifying $b^\si$, we can thus rewrite \eqref{2-3} as
%
%
\begin{eqnarray}
\label{PFE-1} \si_t&=&\si_0+\int_0^tb^\si_s
\,ds+\int_0^tH^\si_s
\,dW_s+\int_0^tH^{\prime \si}_s
\,dW_s'
\nonumber
\\[-8pt]
\\[-8pt]
&&{}+\int_0^t
\int_\R\delta^\si(s,z) (\um-\un) (ds,dz).\nonumber
\end{eqnarray}
It\^o's formula gives us
%
%
\begin{eqnarray}
\label{PFE-2} %
c_t&=&c_0+\int
_0^tb^c_s \,ds+\int
_0^tH^c_s
\,dW_s+\int_0^tH^{\prime c}_s
\,dW_s'\nonumber
\\
&&{}+ \int_0^t
\int_\R\delta^c(s,z) (\um-\un) (ds,dz),
\\
\eqntext{\mbox{where } \cases{ %
b_t^c=2
\si_tb_t^\si+ \bigl(H^\si_t
\bigr)^2+ \bigl(H^{\prime \si}_t \bigr)^2+
\displaystyle\int_\R \delta^\si
(t,z)^2 \,dz,\vspace*{3pt}
\cr
H^{c}_t=2 \si_t
H^\si_t,\qquad H^{\prime {c}}_t=2
\si_t H^{\prime \si}_t,\vspace*{3pt}
\cr
\delta^{c}(t,z)=2
\si_{t-}\delta^\si(t,z)+ \delta^c(t,z)^2,
}}
\end{eqnarray}
and we can thus strengthen and complement \eqref{2-6} as follows:
%
%
\begin{eqnarray}
\label{PFE-4} %
&\bigl |\delta(t,z)\bigr |^r\leq J(z),
\qquad\bigl |\delta^\si(t,z)\bigr |^2\leq J(z) ,\qquad\bigl |
\delta^c(t,z)\bigr |^2\leq J(z),&
\\
&|X_t|+|\si_t|+c_t+|b_t|+\bigl |b^\si_t\bigr |+\bigl |H^\si_t\bigr |+\bigl |H^{\prime \si}_t\bigr |\hspace*{28pt}&\nonumber
\\
&\qquad{}+
\bigl |b^{c}_t\bigr |+\bigl |H^{c}_t\bigr |+\bigl |H^{\prime {c}}_t\bigr |+\bigl |
\gamma^{\pm}_t\bigr |+\bigl |\De Y^{\pm}_t\bigr |\leq K,&
\nonumber
\\
&V=X,c,\si,b,\gamma^{+},\gamma^{-},H^\si,H^c\hspace*{120pt}&\nonumber
\\
&\qquad
\quad\Rightarrow \quad \bigl |\E(V_{t+s}-V_{t}\mid
\f_t)\bigr |+\E \bigl(|V_{t+s}-V_{t}|^2\mid
\f_t \bigr)\leq Ks.&
\nonumber
\end{eqnarray}
Here, $K$ is a constant, and below $K$ and $\phi_n$ will denote a
constant and
a sequence of (nonrandom) numbers going to $0$ as $n\to\infty$, all these
changing from line to line. They may depend on the characteristics of
$X$ and
on the powers for which the forthcoming estimates are stated. Moreover,
in the
theorem to be proven, the arguments $u$ in $\wC(u)^n_t$ or $\wC'(u)^n_t$
are $u=\theta u_n\to0$, where $\theta$ varies in a \emph{fixed} set
$\Te\subset
(0,\infty)$: hence in the sequel we \emph{implicitly assume} $u\in(0,1]$.

Upon replacing $g(x)$ by $g(1)+1$ when $x>1$, we
get \eqref{2-5} for all $x\in(0,\infty)$. We lose the fact that $g$ is
decreasing, but it is still decreasing on $(0,1]$, hence $x^{r-1}g(x)$ as
well because $r\leq1$, and the property
$\int_0^1x^{r-1}g(x) \,dx<\infty$ implies $x^rg(x)\to0$ as $x\to0$.
Summarizing
and recalling $\be\geq1$, we have
%
%
\begin{eqnarray}
\label{PFE-5}  x>0 \quad\Rightarrow\quad \biggl|\BF^{\pm}(x)-
\frac{1}{x^{\be}} \biggr|&\leq& g(x),\qquad \BF^{\pm}(x)\leq\frac{K}{x^\be},\quad\mbox{and}
\nonumber
\\[-8pt]
\\[-8pt]
 \lim_{x\to
0} x^rg(x)&=&0.\nonumber
\end{eqnarray}

Below we unify the proofs of the claims (a) and (b). This is at the expense
of somewhat cumbersome notation, but it saves a lot of space because the
proofs are totally similar. To this end, we introduce a number $\ka$ which
takes the value $1$ if we deal with the nonsymmetrized version and the value
$2$ when we consider the symmetrized version. We set\vspace*{-1pt}
%
%
\begin{eqnarray}
\label{PFE-6} %
L(1,u)^n_j&=&L'(u)^n_j,
\qquad L(2,u)^n_j=L(u)^n_j,\nonumber
\\
 \widehat{c}(1,u)^n_j&=&\widehat{c}'(u)^n_j,
\qquad\widehat {c}(2,u)^n_j=\widehat{c}(u)^n_j,
\\
\wC(\ka,u)^n_t&=&\ka v_n\sum
_{j=0}^{[t/\ka v_n]-1} \biggl(\widehat{c}(\ka,u)^n_j-
\frac{2}{\ka u^2k_n} \bigl(\sinh \bigl(\ka u^2\widehat{c}(
\ka,u)^n_j/2 \bigr) \bigr)^2 \biggr)
\nonumber
\end{eqnarray}
[so $\wC(1,u)^n_t=\wC'(u)^n_t$ and $\wC(2,u)^n_t=\wC(u)^n_t$], and also
(recall that when $\ka=1$ we suppose $\be>1$, so the quantities below
are well defined)\vspace*{-1pt}
%
%
\begin{eqnarray}
\label{PFE-6A} %
A(1,u)_t^n&=&A'(u)^n_t,
\qquad A(2,u)_t^n=A(u)^n_t,
\nonumber
\\[-8pt]
\\[-8pt]
Z(\ka,u)^n_t&=&\frac{1}{\rdn} \bigl(\wC(
\ka,u)^n_t-C_t-A(\ka ,u)^n_t
\bigr).
\nonumber
\end{eqnarray}

Next, recalling the notation \eqref{F-5}, we set\vspace*{-1pt}
%
%
\begin{eqnarray}
\label{PFE-7} %
U(\ka,u)_t&=&e^{-\ka u^2c_t/2},
\qquad \BU(\ka,u)^n_t=e^{-\ka\De_n^{1-\be/2}u^{\be} a_t},
\nonumber
\\
\wU(1,u)^n_t&=&\cos \bigl(\De_n^{1-\be/2}u^{\be}
a_t' \bigr),\qquad \wU(2,u)^n_t=1,
\\
\ua(\ka,u)^n_t&=&U(\ka,u)_t \BU(
\ka,u)^n_t \wU(\ka ,u)^n_t.
\nonumber
\end{eqnarray}
Since $0\leq c_t\leq K$ and $0\leq a_t\leq K$ and $|a_t'|\leq K$,
and assuming $n$ large enough to have $\De_n^{1-\be}u^{\be}|a_t'|
\leq\frac{1}2$ for all $t$ and $u\in(0,1]$, we see that, for some
$\chi\in(0,1)$,\vspace*{-1pt}
%
%
\begin{eqnarray}
\label{PFE-701} \chi&\leq& U(\ka,u)_t\leq1,\qquad\chi\leq\BU(
\ka,u)^n_t\leq 1,
\nonumber
\\[-8pt]
\\[-8pt]
 \chi&\leq&\wU(\ka,u)^n_t
\leq1,\qquad\chi\leq\ua(\ka,u)^n_t\leq1.\nonumber
\end{eqnarray}
Moreover, It\^o's formula yields\vspace*{-1pt}
%
%
\begin{eqnarray}
\label{PFE-3} %
U(\ka,u)_t &=&U(
\ka,u)_0+\int_0^tb^{U(\ka,u)}_s
\,ds+\int_0^tH^{U(\ka,u)}_s
\,dW_s +\int_0^tH^{\prime {U(\ka,u)}}_s
\,dW_s'
\nonumber
\\
&&{}+\int_0^t\int_\R
\delta^{U(\ka,u)}(s,z) (\um-\un) (ds,dz),
\\
&&\mbox{where } \cases{ %
b^{U(\ka,u)}_t\!=-
\displaystyle\frac{\ka u^2}2 U(\ka,u)_t b_t^c +
\frac{\ka^2u^4}8 U(\ka,u)_t \bigl( \bigl(H_t^c
\bigr)^2\!+ \bigl(H_t^{\prime c} \bigr)^2
\bigr)\vspace*{3pt}
\cr
\phantom{b^{U(\ka,u)}_t=} {}+U( \ka,u)_t\!
\displaystyle\int_\R \biggl(\!e^{-\ka
u^2\delta^c(t,z)/2}\!-1+ \frac{\ka u^2}2
\delta^c(t,z)\! \biggr)\,dz,\vspace*{3pt}
\cr
H^{U(\ka,u)}_t=-
\displaystyle\frac{\ka u^2}2 U(\ka,u)_t H^c_t,\vspace*{3pt}
\cr
H^{\prime {U(\ka,u)}}_t=- \displaystyle\frac{\ka u^2}2 U(\ka,u)_t
H^{\prime c}_t,\vspace*{3pt}
\cr
\delta^{U(\ka,u)}(t,z)=U(
\ka,u)_{t-} \bigl(e^{-\ka u^2\delta
^c(t,z)/2}-1 \bigr). 
}
\nonumber
\end{eqnarray}
Therefore, we have for all $u\in(0,1]$:
%
%
\begin{eqnarray}
\label{PFE-8} \bigl |b^{U(\ka,u)}_t\bigr |+\bigl |H^{U(\ka,u)}_t\bigr |+\bigl |H^{\prime {U(\ka,u)}}_t\bigr |
&\leq& Ku^2,
\nonumber
\\[-8pt]
\\[-8pt]
 \bigl |\delta^{U(\ka,u)}(t,z)\bigr |^2&\leq&
Ku^2J(z).\nonumber
\end{eqnarray}

Since $ ||x|^\be-|y|^\be-\be\{y\}^{\be-1}(x-y) |\leq
K|x-y|^\be$
for $x,y\in\R$ when $1\leq\be<2$, and a similar estimate for $\{x\}
^\be
-\{y\}^\be$, and since $|\gamma^{\pm}_t|\leq K$, the last
part of \eqref{PFE-4} implies for $s\in[0,1]$ and $q\geq2$:
%
%
\begin{eqnarray}
\label{PFE-9} %
\bigl |\E(a_{t+s}-a_t\mid
\f_t)\bigr |+\bigl |\E \bigl(a_{t+s}'-a_t'
\mid\f_t \bigr)\bigr | &\leq& Ks^{\be/2},
\nonumber
\\[-8pt]
\\[-8pt]
\E \bigl(|a_{t+s}-a_t|^q+\bigl |a_{t+s}'-a_t'\bigr |^q
\mid\f_t \bigr)&\leq& Ks^{1\wedge(q\be/2)}.
\nonumber
\end{eqnarray}
Using $|e^x-e^y-e^y(x-y)|\leq(x-y)^2$ for $x,y\leq0$ and a similar estimate
for the cosine function when $\ka=1$, we deduce for all $u>0$ and
$q\geq2$:
%
%
\begin{eqnarray}
\label{PFE-10} %
&\bigl |\E \bigl(U(\ka,u)_{t+s}-U(
\ka,u)_t\mid\f_t \bigr)\bigr |\leq Ku^2s,&\nonumber
\\
 &\E
\bigl(\bigl |U(\ka,u)_{t+s}-U(\ka,u)_t\bigr |^q\mid
\f_t \bigr)\leq Ku^{2q}s,&
\nonumber
\\
&\bigl |\E \bigl(\BU(\ka,u)_{t+s}^n-\BU(\ka,u)_t^n
\mid\f_t \bigr)\bigr | +\bigl |\E \bigl(\wU(\ka,u)_{t+s}^n-
\wU(\ka,u)_t^n\mid\f_t \bigr)\bigr |&
\nonumber
\\[-8pt]
\\[-8pt]
&\qquad\leq K
\De_n^{1-\be/2}u^{\be} s^{\be/2},\hspace*{175pt}&\nonumber
\\
&\E \bigl(\bigl |\BU(\ka,u)_{t+s}^n-\BU(\ka,u)_t^n\bigr |^q
+\bigl |\wU(\ka,u)_{t+s}^n-\wU(\ka,u)_t^n\bigr |^q
\mid\f_t \bigr)& \nonumber
\\
&\qquad\leq K\De_n^{q(1-\be/2)}u^{q\be}s.\hspace*{153pt}&
\nonumber
\end{eqnarray}
In turn, since $xy-zw=(x-z)(y-w)+z(y-w)+w(x-z)$, this yields
%
%
\begin{eqnarray}
\label{PFE-11} %
&\quad\bigl |\E \bigl(\ua(\ka,u)^n_{t+s}-
\ua(\ka,u)^n_t\mid\f_t \bigr)\bigr | \leq K
\bigl(u^2s+\De_n^{1-\be/2}u^{\be}
s^{\be/2} \bigr),&
\nonumber
\\
&\quad\E \bigl(\bigl |\ua(\ka,u)^n_{t+s}-\ua(\ka,u)^n_t\bigr |^q
\mid\f_t \bigr)\leq K s \bigl(u^{2q}+\De_n^{q(1-\be/2)}u^{q\be}
\bigr),&
\nonumber
\\[-10pt]
\\[-6pt]
&\quad\E \bigl(\bigl |\ua(\ka,u)^n_{t+s}-\ua(\ka,u)^n_t-
\bigl(U(\ka,u)_{t+s}-U(\ka,u)_t \bigr) \BU(
\ka,u)_t^n \wU(\ka,u)^n_t\bigr |^q
\mid\f_t \bigr)&\nonumber
\\
&\quad\qquad\leq K s \De_n^{q(1-\be/2)}u^{q\be}.\hspace*{224pt}&
\nonumber
\end{eqnarray}

We end this preliminary subsection with another set of notation, with again
$\ka=1,2$.
\begin{eqnarray*}
\rho'(1)^n_i&=&\frac{1}{\rdn}
\si_{(i-1)\De_n}\De^n_iW,\qquad \rho'(2)^n_i=
\frac{1}{\rdn}\si_{(i-1)\De_n} \bigl(\De^n_iW-
\De^n_{i+1}W \bigr),
\\
\rho''(1)^n_i&=&
\frac{1}{\rdn} \bigl(\gamma^{+}_{(i-1)\De_n}\De^n_iY^{+}
+\gamma^{-}_{(i-1)\De_n}\De^n_iY^{-}
\bigr),
\\
\rho''(2)^n_i&=&
\frac{1}{\rdn} \bigl(\gamma^{+}_{(i-1)\De_n} \bigl(\De
^n_iY^{+}-\De^n_{i+1}Y^{+}
\bigr) +\gamma^{-}_{(i-1)\De_n} \bigl(\De^n_iY^{-}-
\De^n_{i+1}Y^{-} \bigr) \bigr),
\\
\rho(\ka)^n_i&=&\rho'(
\ka)^n_i+\rho''(
\ka)^n_i,\qquad \Brho(1)^n_i=
\frac{1}{\rdn}\De^n_iX,
\\
\Brho(2)^n_i&=&
\frac
{1}{\rdn} \bigl(\De^n_iX-\De^n_{i+1}X
\bigr),
\end{eqnarray*}
\begin{eqnarray}
\label{PFE-13}\quad  \xi(\ka,u)^{w,n}_j&=&\cases{ %
\displaystyle\frac{1}{k_n}\sum_{l=0}^{k_n-1}
\bigl(\cos \bigl(u\rho(\ka)^n_{1+\ka
jk_n+\ka l} \bigr) -\ua(
\ka,u)^n_{\ka(jk_n+l)\De_n} \bigr)\vspace*{3pt}
\cr
\qquad\mbox{if } w=1,\vspace*{3pt}
\cr
\displaystyle\frac{1}{k_n}\sum_{l=0}^{k_n-1} \bigl(\cos
\bigl(u\Brho(\ka)^n_{1+\ka
jk_n+\ka l} \bigr) -\cos \bigl(u\rho(
\ka)^n_{1+\ka jk_n+\ka l} \bigr) \bigr)\vspace*{3pt}
\cr
 \qquad\mbox{if } w=2,\vspace*{3pt}
\cr
\displaystyle\frac{1}{k_n}\sum_{l=0}^{k_n-1} \bigl(\ua(
\ka,u)^n_{\ka(jk_n+l)\De_n} -\ua(\ka,u)^n_{\ka jv_n}
\bigr)\vspace*{3pt}
\cr
 \qquad\mbox{if } w=3, 
}
\nonumber
\\
\quad \xi(\ka,u)^{n}_j&=&\frac{1}{\ua(\ka,u)^n_{\ka jv_n}}\sum
_{w=1}^3 \xi(\ka,u)^{w,n}_j,
\\
\quad \Om(\ka,u)_{n,t}&=& \biggl\{\sup_{j=0,\ldots,[t/\ka v_n]-1} \bigl |\xi(
\ka,u)^n_j\bigr | \leq\frac{1}2 \biggr\}.
\nonumber
\end{eqnarray}
Note that, by virtue of \eqref{PFE-701},
%
%
\begin{equation}
\label{PFE-131} \bigl |\xi(\ka,u)^{w,n}_j\bigr |\leq K,\qquad\bigl |\xi(
\ka,u)^{n}_j\bigr |\leq K.
\end{equation}
Finally, let us mention that below we assume \eqref{F-8}. This implies the
following properties, which will be used many times below for various values
of the reals $w_j$ below:
%
%
\begin{equation}
\label{P} \frac{k_n^{w_1} v_n^{w_2} \De_n^{w_3}}{u_n^{w_4}}\to0 \qquad\mbox {if } \cases{ %
\mbox{either}\quad w_3>\displaystyle\frac{w_1-w_2}2
\cr
\mbox{or}
\quad w_4<4(w_1+w_2),\qquad w_3
\geq\displaystyle\frac{w_1-w_2}2. 
}\hspace*{-35pt}
\end{equation}

\subsection{The scheme of the proof}

We have the sequence $u_n$ and $\be\in[1,2)$ with further $\be>1$
when we
deal with (b) of Theorem~\ref{TF-1}, hence when $\ka=1$. Below,
$\theta$
always belongs to a finite set $\Te\subset(0,\infty)$ which,
without loss of generality, contains $1$. We set\vspace*{-1pt}
\begin{eqnarray*}
\Om(\ka)_{n,t}&=&\bigcap_{\theta\in\Te} \Om(\ka,
\theta u_n)_{n,t},\qquad a^{n}_t=
\De_n^{1-\be/2} a_t,
\\[-2pt]
 a^{\prime {n}}_t&=&
\cases{ %
\De_n^{1-\be/2}
a_t'&\quad$\mbox{if }\be>1$,
\cr
0&\quad$\mbox{if }
\be=1$, 
}
\\[-2pt]
f_{\ka,u}(x)&=& \bigl(\sinh \bigl(\ka u^2x/2 \bigr)
\bigr)^2,\qquad h_{1,u} \bigl(x,x' \bigr)=
\frac{2}{u^2} \bigl(u^{\be}x-\log \bigl(\cos \bigl(u^{\be
}x'
\bigr) \bigr) \bigr),
\\[-2pt]
 h_{2,u} \bigl(x,x'
\bigr)&=&2u^{\be-2}x.
\end{eqnarray*}
Because $c_t,a_t,a_t'$ are bounded, we have the estimates (with
$f'$ and $f''$ the first two derivatives of $f$):\vspace*{-2pt}
%
%
\begin{eqnarray}
\label{PFE-15} %
&f_{\ka,u}(c_t)+\bigl |f_{\ka,u}'(c_t)\bigr |+\bigl |f_{\ka,u}''(c_t)\bigr |
\leq Ku^4,&
\nonumber
\\[-1pt]
&\bigl |u^2x\bigr |\leq K \quad\Rightarrow\quad\bigl |u^2f_{\ka,u}'(x)\bigr |+\bigl |f_{\ka
,u}''(x)\bigr |
\leq Ku^4,&
\\[-1pt]
&\De_n\leq Ku^2 \quad\Rightarrow\quad \bigl |h_{\ka,u}
\bigl(a^{n}_t,a^{\prime {n}}_t \bigr)\bigr |\leq
Ku^{\be-2}\De_n^{1-\be/2}\leq K&
\nonumber
\end{eqnarray}
and also\vspace*{-2pt}
%
%
\begin{eqnarray}
\label{PFE-16} -\frac{2}{\ka u^2} \log\ua(\ka,u)_t^n&=&c_t
+ h_{\ka,u} \bigl(a^{n}_t,a^{\prime {n}}_t
\bigr),
\nonumber
\\[-8pt]
\\[-8pt]
 A(\ka,u)^n_t&=&\int_0^th_{\ka,u}
\bigl(a_s^{n},a^{\prime {n}}_s \bigr)
\,ds.\nonumber
\end{eqnarray}

(1) The key step of the proof is as follows. By construction, we have
$L(\ka,u)^n_j=\ua(\ka,u)^n_{jv_n}(1+\xi(\ka,u)^n_j)$. Moreover,
we have $\ua(\ka,\theta u_n)^n_t\geq\chi>0$ by\vspace*{-1pt} \eqref{PFE-701} and
there is
a nonrandom integer $n_0$ such that $k_n\geq4/\chi^2$ for $n\geq n_0$,
implying $L(\ka,\theta u_n)^n_j\geq1/\sqrt{k_n}$ for all $j\leq
[t/v_n]-1$ such
that $1+\xi(\ka,\theta u_n)^n_{j}\geq\frac{1}2$. Hence, we
deduce from \eqref{PFE-16} that\vspace*{-2pt}
%
%
\begin{eqnarray}
\label{PFE-17} %
&&n\geq n_0,\qquad\omega\in
\Om(\ka)_{n,t}
\nonumber
\\[-1pt]
&&\phantom{n\geq n_0,\qquad}\qquad\Rightarrow\quad \widehat{c}(\ka,\theta
u_n)^n_j=c_{\ka jv_n}
+h_{\ka,\theta u_n} \bigl(a^{n}_{\ka jv_n},a^{\prime {n}}_{\ka jv_n}
\bigr)
\\[-1pt]
&&\phantom{n\geq n_0,\qquad\qquad\Rightarrow\quad \widehat{c}(\ka,\theta
u_n)^n_j=}{}- \frac{2}{\ka(\theta u_n)^2} \log \bigl(1+\xi(\ka,\theta u_n)^n_j
\bigr).\nonumber 
\end{eqnarray}

Another key point is as such: on the set $\Om_{n,t}$ and again for
$n\geq n_0$, we can expand $\log(1+x)$ around $0$ and $f_{\ka,u}$ around
$c_{\ka j\De_n}$ to obtain\vspace*{-1pt}
\begin{eqnarray*}
&& \biggl |\widehat{c}(\ka,\theta u_n)^n_j-c_{\ka jv_n}
-h_{\ka,\theta u_n} \bigl(a^{n}_{\ka jv_n},a^{\prime {n}}_{\ka jv_n}
\bigr)
\\
&&\quad{}+\frac{2}{\ka(\theta u_n)^2} \xi(\ka,\theta u_n)^n_j
-\frac{1}{\ka(\theta u_n)^2} \bigl |\xi(\ka,\theta u_n)^n_j\bigr |^2
\biggr |
\\
&&\qquad\leq\frac{K}{u_n^2} \bigl |\xi(\ka,\theta u_n)^n_j\bigr |^3,
\\
&& \biggl |f_{\ka,\theta u_n} \bigl(\widehat{c}(\ka,\theta u_n)^n_j
\bigr)-f_{\theta
u_n}(c_{\ka jv_n})+ \frac{2}{\ka(\theta u_n)^2}
f_{\ka,\theta u_n}'(c_{\ka jv_n})\xi(\ka ,\theta
u_n)^n_j \biggr |
\\
&&\qquad\leq K \bigl(u_n^{\be}\De_n^{1-\be/2}+\bigl |
\xi(\ka,\theta u_n)^n_j\bigr |^2
\bigr),
\end{eqnarray*}
where for the last estimate we have used \eqref{PFE-15} and the
fact that $|\widehat{c}(\ka,\theta u_n)^n_j|\leq K/u_n^2$ [by the
first estimate, plus
again \eqref{PFE-15} and \eqref{PFE-131}], hence
$|u_n^2f_{\ka,\theta u_n}'(x)|+|f_{\ka,\theta u_n}''(x)|\leq Ku_n^4$
for all $x$
between
$\widehat{c}(\ka,\theta u_n)^n_j$ and $c_{\ka jv_n}$. In turn, this and
\eqref{P} yield on the set $\Om(\ka)_{n,t}$ and for $n\geq n_0$ again:
%
%
\begin{eqnarray}
\label{PFE-18} %
&& \biggl| \biggl(\widehat{c}(\ka,\theta
u_n)^n_j-c_{\ka jv_n}-
h_{\ka,\theta u_n} \bigl(a^{n}_{\ka jv_n},a^{\prime {n}}_{\ka jv_n}
\bigr) -\frac{2}{\ka k_n (\theta u_n)^2} f_{\ka,\theta u_n} \bigl(\widehat{c}(\ka ,\theta
u_n)^n_j \bigr) \biggr)
\nonumber\hspace*{-20pt}
\\
&&\quad{}-\frac{2}{\ka(\theta u_n)^2} \biggl(\frac{2}{\ka k_n(\theta u_n)^2} f_{\ka,\theta u_n}'
(c_{\ka jv_n})-1 \biggr)\xi(\ka,\theta u_n)^n_j
\nonumber
\\[-8pt]
\\[-8pt]
&&\hspace*{142pt}{}-\frac{1}{\ka(\theta u_n)^2} \biggl(\bigl |\xi(\ka,\theta u_n)^n_j\bigr |^2
-\frac{2}{k_n} f_{\ka,\theta u_n}(c_{\ka jv_n}) \biggr)\biggr |
\nonumber\hspace*{-20pt}
\\
&&\qquad\leq K \biggl(\frac{|\xi(\ka,\theta u_n)^n_j|^2}{k_nu_n^2}+ \frac{|\xi(\ka,\theta u_n)^n_j|^3}{u_n^2} +\frac{\De_n^{1-\be/2}}{k_nu_n^{2-\be}}
\biggr).
\nonumber\hspace*{-20pt}
\end{eqnarray}

(2) Recalling \eqref{PFE-6A} and \eqref{PFE-16}, we can write
%
%
\begin{equation}
\label{PFE-19} Z(\ka,\theta u_n)^n=V^{\ka,n,\theta}+V^{\prime\ka,n,\theta}+V^{\prime\prime \ka
,n,\theta},
\end{equation}
where
\begin{eqnarray*}
V^{\ka,n,\theta}_t&=&-\frac{1}{\rdn}\int_{\ka v_n([t/\ka
v_n]-1)}^t
\bigl(c_s +h_{\ka,\theta u_n} \bigl(a^{n}_s,a^{\prime {n}}_s
\bigr) \bigr)\, ds,
\\
V^{\prime {\ka,n,\theta}}_t&=&-\sum_{j=0}^{[t/\ka v_n]-1}
\frac{1}{\rdn}\int_{\ka jv_n}^{\ka(j+1)v_n}
\bigl(c_s-c_{\ka jv_n}%
\\
&&\hspace*{124pt}{}+ \bigl(
h_{\ka,\theta u_n} \bigl(a^{n}_s,a^{\prime {n}}_s
\bigr)
\\
&&\hspace*{140pt}{}-h_{\ka,\theta u_n} \bigl(a^{n}_{\ka jv_n},a^{\prime {n}}_{\ka,jv_n}
\bigr) \bigr) \bigr) \,ds,
\\
V^{\prime\prime {\ka,n,\theta}}_t&=&\sum_{j=0}^{[t/\ka v_n]-1}
\frac{\ka
v_n}{\rdn} \biggl(\widehat{c}(\ka, \theta u_n)^n_j-c_{\ka jv_n}%
 h_{\ka,\theta u_n} \bigl(a^{n}_{\ka jv_n},a^{\prime {n}}_{\ka jv_n}
\bigr)
\\
&&\hspace*{94pt}{}-\frac{2}{\ka k_n(\theta u_n)^2} f_{\ka,\theta u_n} \bigl(\widehat{c}(\ka ,\theta
u_n)^n_j \bigr) \biggr).
\end{eqnarray*}

Let us also introduce the following processes:
\begin{eqnarray*}
\BV^{\ka,n,\theta}_t&=&\sum_{j=0}^{[t/\ka v_n]-1}
\frac{2 v_n}{(\theta u_n)^2\rdn} \biggl(\frac{2}{\ka k_n(\theta u_n)^2} f_{\ka,\theta u_n}'
(c_{\ka jv_n})-1 \biggr)\xi(\ka,\theta u_n)^n_j,
\\
\BV^{\prime {\ka,n,\theta}}_t&=&\sum_{j=0}^{[t/\ka v_n]-1}
\frac{v_n}{(\theta u_n)^2\rdn} \biggl(\bigl| \xi(\ka,\theta u_n)^n_j\bigr |^2
-\frac{2}{k_n} f_{\ka,\theta u_n}(c_{\ka jv_n}) \biggr),
\\
R^{\ka,n,\theta}_t&=&\sum_{j=0}^{[t/\ka v_n]-1}
\frac{v_n}{\rdn} \biggl(\frac{|\xi(\ka,\theta u_n)^n_j|^2}{k_nu_n^2}+ \frac{|\xi(\ka,\theta u_n)^n_j|^3}{u_n^2} +
\frac{\De_n^{1-\be/2}}{k_nu_n^{2-\be}} \biggr).
\end{eqnarray*}
By virtue of \eqref{PFE-18}, we then obtain
\begin{eqnarray}
 \bigl | Z(\ka,\theta u_n)^{n,l}_s-V^{\ka,n,\theta}_s-V^{\prime\ka
,n,\theta}_s
-\BV^{\ka,n,\theta}_s-\BV^{\prime {\ka,n,\theta}}_s \bigr |\leq
KR^{\ka
,n,\theta}_t \nonumber
\\
\eqntext{\mbox{on $\Om(\ka)_{n,t}$, for
all $s\leq t$}.}
\end{eqnarray}

Therefore, Theorem~\ref{TF-1} follows from the next four lemmas, where
$Z$ and $\BZ$ are as in \eqref{F-12bA}:

%
\begin{lem}\label{L20}
We have $\PP((\Om(\ka)_{n,t})^c)\to0$.
\end{lem}

%
\begin{lem}\label{L21} We have $\frac{1}{u_n^2} V^{\ka,n,\theta
}\stackrel{{\mathit{u.c.p.}}}{\Longrightarrow}0$ and
$\frac{1}{u_n^2} V^{\prime\ka,n,\theta}\stackrel{{\mathit{u.c.p.}}}{\Longrightarrow}0$.
\end{lem}

%
\begin{lem}\label{L22} We have $\frac{1}{u_n^2} R^{\ka,n,\theta
}_t\stackrel{{\mathit{u.c.p.}}}{\Longrightarrow}0$ and
$\frac{1}{u_n^2} \BV^{\prime\ka,n,\theta}\stackrel{{\mathit{u.c.p.}}}{\Longrightarrow}0$.
\end{lem}

%
\begin{lem}\label{L25}
The processes $ (\BV^{\ka,n,1},(\frac{1}{u_n^2}
(\BV^{\ka,n,\theta}-\BV^{\ka,n,1}))_{\theta\in\Te} )$
converge stably in
law to the limit $ (\ka^{1/2}Z,(\ka^{3/2}(\theta^2-1)\BZ
)_{\theta\in\Te} )$,
provided $\be>1$ when $\ka=1$.
\end{lem}

\subsection{Proofs of Lemmas \texorpdfstring{\protect\ref{L20}--\protect\ref{L25}}{7--10}}

We begin with Lemma~\ref{L21}, which is simple to prove.

\begin{pf*}{Proof of Lemma~\ref{L21}}
By the boundedness of $c_t$
and the
property $\De_n\leq Ku_n^2$, we deduce from \eqref{PFE-15} that
$|V^{\ka,n,\theta}_t|\leq K\frac{v_n}{\rdn}$, which\vspace*{-4pt} is o$(u_n^2)$
by \eqref{F-8},
hence the first\vspace*{-1pt} claim. Next, we have $h_{2,\theta u_n}(a^{n}_s,a^{\prime {n}}_s)
-h_{2,\theta u_n}(a_w,a_w')=\break \frac{2\De_n^{1-\be/2}} {
(\theta u_n)^{2-\be}} (a_s-a_w)$ and also, as soon as
$(\theta u_n)^{\be}|a^{\prime {n}}_w|\leq\frac{1}2$ (hence for all $n$ large enough),
\begin{eqnarray*}
&& \biggl |h_{1,\theta u_n} \bigl(a^{n}_s,a^{\prime {n}}_s
\bigr) -h_{1,\theta u_n} \bigl(a_w,a_w'
\bigr)
\\
&&\quad{}- \frac{2\De_n^{1-\be/2}}{(\theta u_n)^{2-\be}} \bigl( a_s-a_w-
\bigl(a_s'-a_w' \bigr)\tan
\bigl(\De_n^{1-\be/2}(\theta_lu_n)^{\be}
a_w' \bigr) \bigr) \biggr |
\\
&&\qquad\leq K\De_n^{2-\be}u_n^{2\be-2}
\bigl |a_s'-a_w'\bigr |^2.
\end{eqnarray*}
Hence,\vspace*{-2pt} \eqref{PFE-4} for $V=c$ and \eqref{PFE-9} imply that the $j$th
summand $\ze^n_j$ in the definition of $V^{\prime {\ka,n,\theta}}_t$
satisfies in all
cases
\begin{eqnarray*}
\bigl |\E \bigl(\ze^n_j\mid\f_{\ka jv_n} \bigr)\bigr |&\leq&
\frac{K}{\rdn} \bigl(v_n^2+v_n^{1+\be/2}
\De_n^{1-\be/2}u_n^{\be-2} \bigr) =\mathrm{o}
\bigl(v_n u_n^2 \bigr),
\\[1pt]
\E \bigl( \bigl(\ze^n_j \bigr)^2\mid
\f_{\ka jv_n} \bigr)&\leq&\frac{Kv_n^3}{\De_n} \bigl(1+\De_n^{2-\be}u_n^{2\be-4}
\bigr)=\mathrm{o} \bigl(v_nu_n^4 \bigr),
\end{eqnarray*}
where the last two estimates follow from \eqref{P}.
Then a classical argument yields the second claim.
\end{pf*}

The other lemmas need quite many preliminary results. Below, to ease
notation we simply write $u_n$ instead of $\theta u_n$.

%
\begin{lem}\label{L1}
Recalling \eqref{PFE-13}, we have for all
$q\geq2$:
%
%
\begin{eqnarray}
\label{PFE-21} %
\bigl |\E \bigl(\cos \bigl(u_n
\Brho( \ka)^n_i \bigr)-\cos \bigl(u_n\rho(
\ka)^n_i \bigr)\mid \f_{(i-1)\De_n} \bigr) \bigr |&\leq&
u_n^4\rdn\phi_n,
\nonumber
\\[-8pt]
\\[-8pt]
\E \bigl(\bigl |\cos \bigl(u_n\Brho(\ka)^n_i
\bigr)-\cos \bigl(u_n\rho(\ka)^n_i
\bigr)\bigr |^q\mid\f _{(i-1)\De_n} \bigr) &\leq& u_n^4
\rdn\phi_n.
\nonumber
\end{eqnarray}
\end{lem}

\begin{pf} (1) We begin the proof with the case $\ka=1$.
Letting
$\BX_t=\int_0^t\int_\R\delta(s,\allowbreak  z) \um(ds,dz)$,
we have $u_n\Brho(1)^n_i=\sum_{k=1}^4\theta(k)^n_i$, where
\begin{eqnarray*}
\theta(1)^n_i&=&u_n\rho(1)^n_i,
\\[1pt]
 \theta(2)^n_i&=&\frac{u_n}{\rdn} \int
_{(i-1)\De_n}^{i\De_n}(\si_s-\si_{(i-1)\De_n})
\,dW_s +u_n\rdn b_{(i-1)\De_n},
\\[1pt]
\theta(3)^n_i&=&\frac{u_n}{\rdn} \int
_{(i-1)\De_n}^{i\De_n}(b_s-b_{(i-1)\De_n})
\,ds,
\\[1pt]
\theta(4)^n_i&=&\frac{u_n}{\rdn} \biggl(
\De^n_i\BX+ \int_{(i-1)\De_n}^{i\De_n}
\bigl(\gamma^{+}_s-\gamma^{+}_{(i-1)\De_n}
\bigr) \,dY^{+}_s
\\[1pt]
&&\hspace*{55pt}{}+ \int_{(i-1)\De_n}^{i\De_n}
\bigl(\gamma^{-}_s-\gamma^{-}_{(i-1)\De_n}
\bigr) \,dY^{-}_s \biggr).
\end{eqnarray*}
We also write $\Bte(k)^n_i=\sum_{m=1}^k\theta(m)^n_i$, so
%
%
\begin{equation}
\label{PFE-22} \cos \bigl(u_n\rho(1)^n_i
\bigr)=\cos \bigl(\Bte(1)^n_i \bigr),\qquad \cos
\bigl(u_n\Brho(1)^n_i \bigr)=\cos \bigl(
\Bte(4)^n_i \bigr).
\end{equation}

(2) In this step, we prove the following estimates, for any $w\geq2$
and $\ep>0$:
%
%
\begin{eqnarray}
\label{PFE-23} %
\E \bigl(\bigl |\theta(k)^n_i\bigr |^w
\mid\f_{(i-1)\De_n} \bigr)&\leq&\cases{ %
Ku_n^w\De_n^{1-w/2}&\quad$\mbox{if }
k=1$,\vspace*{2pt}
\cr
Ku_n^w\De_n&\quad$\mbox{if } k=2$,\vspace*{2pt}
\cr
Ku_n^w\De_n^{1+w/2}&\quad$
\mbox{if } k=3$, 
}
\\
\E \bigl(\bigl |\theta(4)^n_i\bigr |\wedge1\mid\f_{(i-1)\De_n}
\bigr)&\leq& u_n^4\rdn\phi_n,
\nonumber
\\
\E \bigl( \bigl(\bigl |u_n\rho''(
\ka)^n_i\bigr |\wedge1 \bigr)^2\mid
\f_{(i-1)\De_n} \bigr) &\leq& K\De_n^{1-\ep-\be/2}.
\nonumber
\end{eqnarray}
We classically have $\E(|\De^n_iW|^w\mid\f_{(i-1)\De_n})\leq K\De
_n^{w/2}$,
whereas $\E(|\De^n_iY^\pm|^w\mid\f_{(i-1)\De_n})\leq K\De_n$ by
Lemma~2.1.5 of
\cite{JP} (because $Y^\pm$ has bounded jumps), yielding case $k=1$.
Cases $k=2,3$ follow from \eqref{PFE-4}.

For case $k=4$, it is enough to prove the result for each of the three
summands in the definition of $\theta(4)^n_i$. For the first summand
$\De^n_i\BX$, we observe that $|u_n\De^n_i\BX/\rdn|\wedge1\leq
K(|\De^n_i
\BX
/\rdn|\wedge1)$. Then we apply Corollary~2.1.9-(c) of \cite{JP}
with $q=\frac{1}2$ and $s=\De_n$ and $r$ as in \ref{ass(A)} and~\ref{ass(B)}
and $p=1$, to
obtain
%
%
\begin{equation}
\label{PFE-24} \E \biggl(\frac{|u_n\De^n_i\BX|}{\rdn}\wedge1\Bigm|\f_{(i-1)\De
_n}
\biggr)\leq \De_n^{1-r/2}\phi_n.
\end{equation}

The\vspace*{-1pt} other two summands are treated analogously, and we consider only
one of
them, say $\al^n_i=\int_{(i-1)\De_n}^{i\De_n}(\gamma^{+}_s-\gamma
^{+}_{(i-1)\De_n})
\,dY^{+}_s$. We observe that the jump measure of $Y^{+}$,
say $\um'$, is Poisson with compensator $\un'(dt,dz)=\break dt\otimes F^{+}(dz)$,
and $\al^n_i=\De^n_i(\delta*(\um'-\un'))$ if we take $\delta
(t,z)=(\gamma^{+}_{t-}
-\break \gamma^{+}_{(i-1)\De_n})z 1_{\{t>(i-1)\De_n\}}$. The notation (2.1.35)
of \cite{JP} for $\delta*(\um'-\un')$ becomes
\begin{eqnarray*}
\wde(p,a)_{(i-1)\De_n,i\De_n}&=&\frac{1}{\De_n} \int_{(i-1)\De
_n}^{i\De_n}\!
\bigl |\gamma^{+}_s-\gamma^{+}_{(i-1)\De_n}\bigr |^p
\,ds \,\int_0^{a/|\gamma^{+}_s-\gamma^{+}_{(i-1)\De_n}|}\!z^p
F^{+}(dz),
\\
\wde'(p)_{(i-1)\De_n,i\De_n}&=&\wde(p,1)_{(i-1)\De_n,j\De_n}
\\
&&{}+
\frac{1}{\De_n} \int_{(i-1)\De_n}^{i\De_n}\bigl |
\gamma^{+}_s -\gamma^{+}_{(i-1)\De_n}\bigr | \,ds
\, \int_{1/|\gamma^{+}_s-\gamma^{+}_{(i-1)\De_n}|}^{\infty}z F^{+}(dz),
\\
\wde''(p)_{(i-1)\De_n,i\De_n}&=&\wde(p,1)_{(i-1)\De_n,j\De_n}
\!+\frac{1}{\De_n} \int_{(i-1)\De_n}^{i\De_n}
\!\BF^{+} \bigl(1/\bigl |\gamma^{+}_s-
\gamma^{+}_{(i-1)\De_n}\bigr | \bigr) \,ds,
\end{eqnarray*}
and we observe that, since $\gamma^{+}$ is bounded and $F^{+}$ is
supported by
$[0,A]$ for some finite $A$, necessarily $\wde'(p)_{(i-1)\De_n,i\De
_n}\leq
K\wde''(p)_{(i-1)\De_n,i\De_n}$. \eqref{PFE-5} yields
$\int_0^xz^p F^{+}(dz)\leq Kx^{p-\be}$ when $p>\be$, hence
%
%
\begin{eqnarray}
\label{PFE-25} %
\wde(p,a)_{(i-1)\De_n,i\De_n}&\leq&
Ka^{p-\be}\frac{1}{\De_n} \int_{(i-1)\De_n}^{i\De_n}
\bigl |\gamma^{+}_s-\gamma^{+}_{(i-1)\De_n}\bigr |^\be
\,ds,
\nonumber\hspace*{-35pt}
\\[-8pt]
\\[-8pt]
\wde'(p)_{(i-1)\De_n,i\De_n}&\leq& K\wde''(p)_{(i-1)\De_n,i\De
_n}
\leq \frac{K}{\De_n} \int_{(i-1)\De_n}^{i\De_n} \bigl |
\gamma^{+}_s-\gamma^{+}_{(i-1)\De_n}\bigr |^\be
\,ds.
\nonumber\hspace*{-35pt}
\end{eqnarray}

Since $1\leq\be<2$, we then use Lemma~2.1.6 of \cite{JP} with
$q=\frac{1}2$
and $r=p\in(\be,2]$ and $s=\De_n$. Since $\E(|\gamma^{+}_{(i-1)\De_n+s}
-\gamma^{+}_{(i-1)\De_n}|^{\be}\mid\f_{(i-1)\De_n})
\leq Ks^{\be/2}$ by \eqref{PFE-4}, we obtain for $p>\be$:
%
%
\begin{equation}
\label{PFE-26} \E \biggl( \biggl(\frac{|u_n\De^n_i(\delta*(\um'-\un'))|}{\rdn
}\wedge 1
\biggr)^p \Bigm|\f_{(i-1)\De_n} \biggr)\leq K\De_n^{1-(p-\be)/4}.
\end{equation}
We then apply H\"older's inequality to get
\[
\E \biggl(\frac{|u_n\al^n_i|}{\rdn}\wedge1\Bigm|\f_{(i-1)\De
_n} \biggr) \leq K
\De_n^{1/p-1/4+\be/4p}.
\]
Under \eqref{P}, both this and \eqref{PFE-24} are smaller than
$u_n^4\rdn\phi_n$, upon choosing $p$ close enough to $\be$ above. Hence,
\eqref{PFE-23} holds for $k=4$.

Finally, the last estimate in \eqref{PFE-23} is obtained exactly as
above, upon
taking $\gamma^{+}_{(i-1)\De_n}$ instead of
$\gamma^{+}_{t-}-\gamma^{+}_{(i-1)\De_n}$,
so the bounds in \eqref{PFE-25} become $Ka^{p-\be}$ and $K$, and the
one in
\eqref{PFE-26} is $K\De_n^{1-(p+\be)/4}$. We then apply the latter
with $p$ close enough to $\be$, and the result follows.

(3) Since $|\cos(x+y)-\cos(x)|\leq1\wedge|y|\wedge(|xy|+y^2)$ and
$|\cos(x+y)-\cos(x)-y\sin(x)|\leq Ky^2$, we deduce from
\eqref{PFE-23} and the Cauchy--Schwarz inequality that
\begin{eqnarray*}
&w\geq1 \quad\Rightarrow\quad \E \bigl(\bigl |\cos \bigl(\Bte(4)^n_i
\bigr)-\cos \bigl(\Bte(3)^n_i \bigr)\bigr |^w
\mid \f_{(i-1)\De_n} \bigr) \leq u_n^4\rdn
\phi_n,&
\\
&\E \bigl(\bigl |\cos \bigl(\Bte(3)^n_i \bigr)-\cos \bigl(
\Bte(2)^n_i \bigr)\bigr |^2\mid\f_{(i-1)\De
_n}
\bigr)\hspace*{37pt}&
\\
&{}+\E \bigl(\bigl |\cos \bigl(\Bte(2)^n_i \bigr)-\cos
\bigl( \theta(1)^n_i \bigr)\bigr |^2\mid\f
_{(i-1)\De_n} \bigr)&
\\
&\leq Ku_n^2\De_n,\hspace*{119pt}&
\\
&\E \bigl(\bigl |\cos \bigl(\Bte(3)^n_i \bigr)-\cos \bigl(
\Bte(2)^n_i \bigr)\bigr |\mid\f_{(i-1)\De
_n} \bigr) \leq
Ku_n^2\De_n,&
\\
&\E \bigl(\bigl |\cos \bigl(\Bte(2)^n_i \bigr)-\cos \bigl(
\theta(1)^n_i \bigr)-\theta(2)^n_i
\sin \bigl(\theta(1)^n_i \bigr) \bigr |\mid\f_{(i-1)\De_n}
\bigr)\leq Ku_n^2\De_n.&
\end{eqnarray*}
This with $w=2$ and \eqref{PFE-22} and $\rdn=$ o$(u_n^2)$ yield the second
estimate \eqref{PFE-21} for $q=2$, hence for all $q\geq2$ because
$|\cos x|\leq1$, and also (with $w=1$ above) that, for the first estimate,
it only remains to prove that
\[
\bigl |\E \bigl(\theta(2)^n_i \sin \bigl(
\theta(1)^n_i \bigr) \mid\f_{(i-1)\De_n} \bigr)\bigr  |\leq
u_n^4\rdn\phi_n.
\]

Now, we have $|\sin(\theta(1)^n_i)-\sin(u_n\rho'(1)^n_i)|\leq
K(|u_n\rho''(1)^n_i|\wedge1)$, and thus
\begin{eqnarray*}
\E \bigl( \bigl | \theta(2)^n_i \bigl(\sin \bigl(
\theta(1)^n_i \bigr)-\sin \bigl(u_n
\rho'(1)^n_i \bigr)\bigr) \bigr |\mid
\f_{(i-1)\De_n} \bigr)&\leq& Ku_n\De_n^{1-\ep/2-\be/4}
\\
&=&u_n^4\rdn\phi_n
\end{eqnarray*}
by the Cauchy--Schwarz inequality and \eqref{PFE-23}, and where the last equality
comes from \eqref{P}, upon choosing $\ep<\frac{2-\be}4$. Hence, it remains
to prove that
\[
\bigl | \E \bigl(\theta(2)^n_i \sin \bigl(u_n
\rho'(1)^n_i \bigr)\mid\f_{(i-1)\De_n}
\bigr) \bigr |\leq u_n^4\rdn\phi_n.
\]

(4) Recalling \eqref{PFE-2}, we set
\[
V_t=\int_0^tH^{\prime \si}_s
\,dW'_s+ \int_0^t
\int_E\delta^\si(s,z) (\um-\un) (ds,dz).
\]
We have the decomposition $\theta(2)^n_i=-\sum_{j=1}^5\mu(j)^n_i$, where
\begin{eqnarray*}
\mu(1)^n_i&=&u_n\rdn b_{(i-1)\De_n},
\\
\mu(2)^n_i&=&\frac{u_n}{\rdn} \int
_{(i-1)\De_n}^{i\De_n} \biggl(\int_{(i-1)\De_n}^sb^\si_t
\,dt \biggr) \,dW_s,
\\
\mu(3)^n_i&=&\frac{u_n}{\rdn} H^\si_{(i-1)\De_n}
\int_{(i-1)\De_n}^{i\De_n}(W_s-W_{(i-1)\De_n})
\,dW_s,
\\
\mu(4)^n_i&=&\frac{u_n}{\rdn} \int
_{(i-1)\De_n}^{i\De_n} \biggl(\int_{(i-1)\De_n}^s
\bigl(H^\si_t-H^\si_{(i-1)\De_n} \bigr)
\,dW_t \biggr) \,dW_s,
\\
\mu(5)^n_i&=&\frac{u_n}{\rdn} \int
_{(i-1)\De_n}^{i\De_n}(V_s-V_{(i-1)\De_n})
\,dW_s,
\end{eqnarray*}
and it thus suffices to prove that, for $j=1,2,3,4,5$:
%
%
\begin{equation}
\label{FPE-27} \bigl | \E \bigl(\mu(j)^n_i \sin
\bigl(u_n\rho'(1)^n_i \bigr)
\mid\f_{{i-1}\De_n} \bigr) \bigr |\leq Ku_n^4\rdn
\phi_n.
\end{equation}

First, $\E (\mu(j)^n_i \sin (u_n\rho'(1)^n_i )\mid\f
_{(i-1)\De_n}
)=0$ for $j=1,3$ follows from the fact that in these
cases the variable whose conditional expectation is taken is
a function of $(\omega,(W_{(i-1)\De_n+t}-W_{(i-1)\De_n})_{t\geq0})$
which is $\f_{(i-1)\De_n}$-measurable\vspace*{1pt} in $\omega$ and odd in the
second argument.
Second, we have $\E((\mu(j)^n_i)^2\mid\break \f_{(i-1)\De_n})\leq
Ku_n^2\De_n^2$ for
$j=2,4$ [use \eqref{PFE-4}], implying
\eqref{FPE-27} for $j=2,4$ by the Cauchy--Schwarz inequality and \eqref{P}.

For analyzing the case $j=5$, we use the representation theorem for
martingales of the Brownian filtration. This implies that the variable\break \vspace*{-1pt}
$\sin(u_n\rho'(1)^n_i)$, whose $\f_{(i-1)\De_n}$-conditional expectation
vanishes, has the form $\int_{i\De_n}^{(i+1)\De_n}L^n_s \,dW_s$ for\vspace*{-1pt}
some process $L^n$, adapted to the filtration $(\f^W_t)_{t\geq0}$ generated
by the process $W$, hence
\begin{eqnarray*}
&&\E \bigl(\mu(5)^n_i \sin \bigl(u_n
\rho'(1)^n_i \bigr)\mid\f_{(i-1)\De_n}
\bigr)
\\
&&\qquad= \frac{u_n}{\rdn} \int_{(i-1)\De_n}^{i\De_n} \E
\bigl((V_s-V_{(i-1)\De_n}) L^n_s\mid
\f_{(i-1)\De_n} \bigr) \,ds.
\end{eqnarray*}
Since further the martingale $V$ is orthogonal to $W$, and by using
once more
the representation theorem [so $L^n_s=\E(L^n_s\mid\f^n_{i-1})+
\int_{(i-1)\De_n}^sL^{\prime n}_t \,dW_t$ for $s\geq i\De_n$],
we deduce $\E((V_s-V_{(i-1)\De_n}) L^n_s\mid\f_{(i-1)\De_n})=0$,
hence $\E (\mu(5)^n_i \sin(u_n\rho^{\prime n}_i)\mid\f_{(i-1)\De
_n} )=0$ and
\eqref{FPE-27} holds for $j=5$. This completes the proof for the case
$\ka=1$.

(5) When $\ka=2$, we do as above, with a few changes:
First $u_n\Brho(2)^n_i=\sum_{k=1}^4\theta(k)^n_i$, where
\begin{eqnarray*}
\theta(1)^n_i&=&u_n\rho(2)^n_i,
\\
\theta(2)^n_i&=&\frac{u_n}{\rdn} \biggl( \int
_{(i-1)\De_n}^{i\De_n}(\si_s-\si_{(i-1)\De_n})
\,dW_s -\int_{i\De_n}^{(i+1)\De_n}(
\si_s-\si_{(i-1)\De_n}) \,dW_s \biggr),
\\
\theta(3)^n_i&=&\frac{u_n}{\rdn} \biggl( \int
_{(i-1)\De_n}^{i\De_n}(b_s-b_{s+\De_n})
\,ds -\int_{i\De_n}^{(i+1)\De_n}(b_s-b_{s+\De_n})
\,ds \biggr),
\\
\theta(4)^n_i&=&\frac{u_n}{\rdn} \biggl(
\De^n_i\BX-\De^n_{i+1}
\BX%
+\int_{(i-1)\De_n}^{i\De_n} \bigl(
\gamma^{+}_s-\gamma^{+}_{(i-1)\De_n} \bigr)
\,dY^{+}_s
\\
&&\hspace*{31pt}{}-\int_{i\De_n}^{(i+1)\De_n}
\bigl(\gamma^{+}_s-\gamma^{+}_{(i-1)\De_n}
\bigr) \,dY^{+}_s
\\
&&\hspace*{31pt}{}+ \int_{(i-1)\De_n}^{i\De_n} \bigl(\gamma^{-}_s-
\gamma^{-}_{(i-1)\De
_n} \bigr)\, dY^{-}_s
\\
&&\hspace*{96pt}{} -
\int_{i\De_n}^{(i+1)\De_n} \bigl(\gamma^{-}_s-
\gamma^{-}_{(i-1)\De_n} \bigr) \,dY^{-}_s
\biggr).
\end{eqnarray*}
The estimates \eqref{PFE-23} remain trivially valid, as well as Step
3. In
Step 4, we use the decomposition $\theta(2)^n_i=-\sum_{j=2}^5\mu
(j)^n_i$, where
\begin{eqnarray*}
\mu(2)^n_i&=&\frac{u_n}{\rdn} \biggl(\int
_{(i-1)\De_n}^{i\De_n} \!\biggl(\int_{(i-1)\De_n}^s\!b^\si_t
\,dt \biggr) \,dW_s -\int_{i\De_n}^{(i+1)\De_n}\!
\biggl(\int_{(i-1)\De_n}^s\!b^\si_t
\,dt \biggr) \,dW_s \biggr),
\\
\mu(3)^n_i&=&\frac{u_n}{\rdn} H^\si_{(i-1)\De_n}
\biggl( \int_{(i-1)\De_n}^{i\De_n}(W_s-W_{(i-1)\De_n})
\,dW_s
\\
&&\hspace*{70pt}{}-\int_{i\De_n}^{(i+1)\De_n}(W_s-W_{(i-1)\De_n})
\,dW_s \biggr),
\\
\mu(4)^n_i&=&\frac{u_n}{\rdn} \biggl( \int
_{(i-1)\De_n}^{i\De_n} \biggl(\int_{(i-1)\De_n}^s
\bigl(H^\si_t-H^\si_{(i-1)\De_n} \bigr)\,
dW_t \biggr) \,dW_s
\\
&&\hspace*{32pt}{}-\int_{i\De_n}^{(i+1)\De_n} \biggl(\int_{(i-1)\De_n}^s
\bigl(H^\si_t-H^\si_{(i-1)\De_n} \bigr)
\,dW_t \biggr) \,dW_s \biggr),
\\
\mu(5)^n_i&=&\frac{u_n}{\rdn} \biggl( \int
_{(i-1)\De_n}^{i\De_n}(V_s-V_{(i-1)\De_n})\,
dW_s -\int_{i\De_n}^{(i+1)\De_n}(V_s-V_{(i-1)\De_n})
\,dW_s \biggr)
\end{eqnarray*}
[the term $\mu(1)^n_i$ no longer shows up]. The rest of proof carries
over without modification.
\end{pf}

%
\begin{lem}\label{L2} We have for all $q\geq2$, and if $u_n'\asymp u_n$:
%
%
\begin{equation}
\label{PFE-32} \bigl | \E \bigl(\cos \bigl(u_n\rho(\ka)^n_i
\bigr)\mid\f_{(i-1)\De_n} \bigr) -\ua(\ka,u_n)^n_{(i-1)\De_n}
\bigr | \leq\phi_n u_n^4 \rdn,
\end{equation}
%
%
\begin{eqnarray}
\label{PFE-33} %
\qquad&& \bigl | \E \bigl(\cos \bigl(u_n
\rho(\ka)^n_i \bigr) \cos \bigl(u_n'
\rho(\ka)^n_i \bigr)\mid \f _{(i-1)\De_n} \bigr)
\nonumber
\\[-4pt]
\\[-12pt]
&&\qquad{}-\tfrac{1}{2} \bigl(\ua \bigl(\ka,u_n+u_n'
\bigr)^n_{(i-1)\De_n} +\ua \bigl(\ka,\bigl |u_n-u_n'\bigr |
\bigr)^n_{(i-1)\De_n} \bigr) \bigr |\leq\phi_nu_n^4
\rdn ,
\nonumber
\end{eqnarray}
%
%
\begin{equation}
\label{PFE-34} \E \bigl(\bigl |\cos \bigl(u_n\rho(\ka)^n_i
\bigr)-\ua(\ka,u_n)^n_{(i-1)\De_n}\bigr |^q\mid
\f _{(i-1)\De_n} \bigr)\leq Ku_n^4.
\end{equation}
\end{lem}

\begin{pf} (1) The variables $\De^n_iW/\rdn$, $\De^n_iY^{+}/\rdn$ and
$\De^n_iY^{-}/\rdn$ are independent one from another and from $\f
_{(i-1)\De_n}$,
with characteristic
functions $\exp(-u^2/2)$ and $\exp(-G^{\pm}_n(u)-iH_n^{\pm}(u))$, where\vspace*{1pt}
\begin{eqnarray*}
G_n^{\pm}(y)&=&\De_n\int_0^1
\biggl(1-\cos\frac{xy}{\rdn} \biggr) F^{\pm}(dx),
\\[1pt]
H_n^{\pm}(y)&=&\De_n\int_0^1
\biggl(\frac{xy}{\rdn}-\sin\frac
{xy}{\rdn} \biggr) F^{\pm}(dx).
\end{eqnarray*}
Analogously, the characteristic functions of $(\De^n_iW-\De
^n_{i+1}W)/\rdn$ and
$(\De^n_iY^{\pm}-\De^n_{i+1}Y^{\pm})/\rdn$ are $\exp(-u^2)$ and
$\exp(-2G^{\pm}_n(u))$. Therefore, by the definition of $\rho(\ka
)^n_i$, and
since $\si_{(i-1)\De_n}$ and
$\gamma^{\pm}_{(i-1)\De_n}$ are $\f_{(i-1)\De_n}$-measur\-able, we have\vspace*{1pt}
%
%
\begin{eqnarray}
\label{PFE-35} %
\E \bigl(\cos \bigl(u_n
\rho(1)^n_i \bigr)\mid\f_{(i-1)\De_n} \bigr)&=&
U(1,u_n)_{(i-1)\De_n} e^{-G^{+}_n(u_n\gamma^{+}_{(i-1)\De_n})-G^{-}_n(u_n\gamma
^{-}_{(i-1)\De_n})}
\nonumber
\\[1pt]
&&{}\times\cos \bigl(H^{+}_n \bigl(u_n
\gamma^{+}_{(i-1)\De_n} \bigr) +H^{-}_n
\bigl(u_n\gamma^{-}_{(i-1)\De_n} \bigr) \bigr),
\\
\E \bigl(\cos \bigl(u_n\rho(2)^n_i \bigr)
\mid \f_{(i-1)\De_n} \bigr)&=&U(2,u_n)_{(i-1)\De_n}
e^{-2G^{+}_n(u_n\gamma^{+}_{(i-1)\De_n})-2G^{-}_n(u_n\gamma
^{-}_{(i-1)\De_n})}.
\nonumber
\end{eqnarray}

(2) In this step, we analyze the behavior of $G^{\pm}_n(y)$ when $y\in(0,A]$
for some $A>0$. Let $\ze_n=\De_n^\eta$ for some $\eta\in(0,\frac{1}2)$, to be
chosen later, so that $\ze_n\to0$
and $\ze_n'=\ze_ny/\rdn\to\infty$. Using \eqref{PFE-5}, we first
see that\vspace*{1pt}
\[
0\leq\int_{\ze_n}^1 \biggl(1-\cos
\frac{xy}{\rdn} \biggr) F^{\pm}(dx) \leq2\BF^{\pm}(
\ze_n)\leq\frac{K}{\ze_n^{\be}}.
\]
Next, Fubini's theorem and a change of variable yield\vspace*{1pt}
\begin{eqnarray*}
&&\int_0^{\ze_n} \biggl(1-\cos\frac{xy}{\rdn}
\biggr) F^{\pm}(dx)
\\[1pt]
&&\qquad=\int_0^{\ze_n'}
\BF^{\pm} \biggl(\frac{z\rdn}{y} \biggr) \sin(z) \,dz -\int
_0^{\ze_n'}\BF^{\pm}(\ze_n)
\sin(z) \,dz,
\end{eqnarray*}
and the absolute value of the last term above is again smaller than
$K/\ze_n^{\be}$ because $ |\int_0^x\sin z \,dz |\leq2$ for
all $x$.
To evaluate the first term, we use \eqref{PFE-5} again to get\vspace*{1pt}
\begin{eqnarray*}
&&\biggl|\int_0^{\ze_n'}\BF^{\pm} \biggl(
\frac{z\rdn}{y} \biggr) \sin (z) \,dz- \frac{y^{\be}}{\De_n^{\be/2}} \chi(\be) \biggr|
\\[1pt]
&&\qquad\leq \biggl|
\frac{y^{\be}}{\De_n^{\be/2}} \int_{\ze_n'}^\infty
\frac{\sin z}{z^{\be}} \,dz \biggr|+ \int_0^{\ze_n'}g \biggl(
\frac{z\rdn}{y} \biggr) \,dz.
\end{eqnarray*}
We have $\int_x^\infty\frac{\sin z}{z^{\be}} \,dz=
\frac{\cos x}{x^{\be}}-\be\int_x^\infty\frac{\cos z}{z^{1+\be}}
\,dz$ by
integration by parts, yielding $ |\int_x^\infty\frac{\sin
z}{z^{\be}} \,dz
|\leq2/x^{\be}$. We also have
\begin{eqnarray*}
\int_0^{\ze_n'}g \biggl(\frac{z\rdn}{y}
\biggr) \,dz&=& \frac{y}{\rdn}\int_0^{\ze_n}g(z)
\,dz \leq\frac{y}{\rdn} \ze_n^{1-r}\int
_0^{\ze_n}\frac{g(z)}{z^{1-r}} \,dz
\\
&=&
\frac{y}{\rdn} \ze_n^{1-r}\phi_n
\end{eqnarray*}
because $\ze_n\to0$. Putting all these together yields
\[
\bigl | G^{\pm}_n(y)-\De_n^{1-\be/2}
y^{\be}\chi(\be)\bigr  |\leq \frac{K\De_n}{\ze_n^{\be}}+\rdn y \ze_n^{1-r}
\phi_n \leq K\De_n^{1-\eta\be}+y\De_n^{1/2+\eta(1-r)}
\phi_n
\]
for all $y>0$, and also (trivially) when $y=0$. Now, we take
$\eta=\frac{1}{2(1-r+\be)}$ and use \eqref{P} to deduce
\[
\bigl |G^{+}_n \bigl(u_n\gamma^{+}_{(i-1)\De_n}
\bigr)+G^{-}_n \bigl(u_n\gamma
^{-}_{(i-1)\De_n} \bigr) -\De_n^{1-\be/2}u_n^{\be}
a_{(i-1)\De_n} \bigr |\leq u_n^4\rdn\phi_n.
\]

Using once more $|e^{x}-e^{y}|\leq|x-y|\wedge1$ if $x,y\leq0$, and
recalling the definition of $\BU(u)^n_t$, we deduce
%
%
\begin{equation}
\label{PFE-36}\qquad  \bigl | e^{-\ka(G^{+}_n(u_n\gamma^{+}_{(i-1)\De_n})
+G^-_n(u_n\gamma^-_{(i-1)\De_n}))}-\BU(\ka,u_n)^n_{(i-1)\De_n}
\bigr | \leq u_n^4\rdn\phi_n.
\end{equation}

(3) Next, we analyze $H^{\pm}_n(y)$: this is for the case when $\ka=1$,
hence $\be>1$. The following estimates are
easy consequences of \eqref{PFE-5}:
\[
0<z\leq1 \quad\Rightarrow\quad \int_0^zx^3
F^{\pm}(dx)\leq Kz^{3-\be},\qquad \int_z^1x
F^{\pm}(dx)\leq Kz^{1-\be}.
\]
With $\ze_n$ and $\ze_n'$ as in the
previous step, and we have
\[
0\leq\int_{\ze_n}^1 \biggl(\frac{xy}{\rdn}-\sin
\frac{xy}{\rdn
} \biggr) F^{\pm}(dx)\leq\frac{y}{\rdn}
\BF^{\pm}(\ze_n)\leq\frac
{Ky}{\rdn\ze_n^{\be}},
\]
\begin{eqnarray*}
&&\int_0^{\ze_n} \biggl(\frac{xy}{\rdn}-\sin
\frac{xy}{\rdn} \biggr) F^{\pm}(dx)
\\
&&\qquad =\int_0^{\ze_n'}
\BF^{\pm} \biggl(\frac{z\rdn}{y} \biggr) (1-\cos z) \,dz -\int
_0^{\ze_n'}\BF^{\pm}(\ze_n)
(1-\cos z) \,dz,
\end{eqnarray*}
and the absolute value of the last term above is smaller than
$Ky\ze_n'/\ze_n^{\be}$. We also have
\begin{eqnarray*}
&&\biggl|\int_0^{\ze_n'}\BF^{\pm} \biggl(
\frac{z\rdn}{y} \biggr) (1-\cos z) \,dz- \frac{y^{\be}}{\De_n^{\be/2}}
\chi'(\be) \biggr|
\\
&&\qquad\leq\frac{y^{\be}}{\De_n^{\be/2}} \int_{\ze_n'}^\infty
\frac{1-\cos z}{z^{\be}} \,dz+ \int_0^{\ze_n'}g \biggl(
\frac{z\rdn}{y} \biggr)\, dz.
\end{eqnarray*}
As seen before, the last term above is less than $\frac{y}{\rdn} \ze_n^{1-r}
\phi_n$, whereas
$\int_x^\infty(1-\allowbreak  \cos z)/{z^\be} \,dz\leq
K/x^{\be-1}$.
Putting all these together, plus $\ze_n'=y\ze_n/\rdn$,
yields for $y>0$:
\[
\bigl |H^{\pm}_n(y)-\De_n^{1-\be/2}
y^{\be}\chi'(\be) \bigr |\leq \frac{Ky\rdn}{\ze_n^{\be}}\leq Ky
\De_n^{1/2-\be\eta}.
\]
The same holds with $-|y|^{\be}$ and $|y|$ instead of $y^{\be}$ and
$y$ when
$y<0$, and it trivially holds for $y=0$. Since
$|\cos x-\cos y|\leq2|x-y|(|x-y|+|y|)$ for all $x,y$, we obtain
\begin{eqnarray*}
&&\bigl |\cos \bigl(H^{+}_n \bigl(u_n
\gamma^{+}_{(i-1)\De_n} \bigr) +H^{-}_n
\bigl(u_n\gamma^{-}_{(i-1)\De_n} \bigr) \bigr)-
\wU(1,u_n)^n_{(i-1)\De_n} \bigr |
\\
&&\qquad\leq K \bigl(u_n^2
\De_n^{1-2\be\eta} +u_n^{1+\be}
\De_n^{3/2-\be/2-\be\eta} \bigr).
\end{eqnarray*}
In view of \eqref{P}, and upon choosing
$\eta>0$ small enough, we deduce that
%
%
\begin{eqnarray}
\label{PFE-37} &&\bigl |\cos \bigl(H^{+}_n \bigl(u_n
\gamma^{+}_{(i-1)\De_n} \bigr) +H^{-}_n
\bigl(u_n\gamma^{-}_{(i-1)\De_n} \bigr) \bigr)-
\wU(1,u_n)^n_{(i-1)\De_n} \bigr |
\nonumber
\\[-8pt]
\\[-8pt]
\nonumber
&&\qquad\leq u_n^4
\rdn\phi_n.
\end{eqnarray}

(4) At this stage, \eqref{PFE-32} is an easy consequence of \eqref{PFE-7},
\eqref{PFE-35}, \eqref{PFE-36} and~\eqref{PFE-37}. Since
\[
\cos \bigl(u_n\rho(\ka)^n_i \bigr) \cos
\bigl(u_n'\rho(\ka)^n_i
\bigr)= \tfrac{1}{2} \bigl(\cos \bigl( \bigl(u_n+u_n'
\bigr)\rho(\ka)^n_i \bigr) +\cos \bigl(\bigl |u_n-u_n'\bigr |
\rho(\ka)^n_i \bigr) \bigr),
\]
\eqref{PFE-33} follows from \eqref{PFE-32}.

Finally, since $|\cos x|\leq1$ and $|\ua(\ka,u)^n_t|\leq1$, it is
enough to
prove \eqref{PFE-34} for $q=2$. Since
$(\cos x)^2=\frac{1}2 (1+\cos(2x))$, an application of \eqref{PFE-32} and
\eqref{PFE-33} shows that the left-hand side of \eqref{PFE-34} is, up
to a remainder term of size smaller than $\phi_nu_n^4\rdn$, equal to
\[
\tfrac{1}{2} \bigl(\ua(\ka,2u_n)^n_{(i-1)\De_n}
-2 \bigl(\ua(\ka,u_n)^n_{(i-1)\De_n}
\bigr)^2+1 \bigr).
\]
An expansion near $0$ of the function $u\mapsto\ua(\ka,u)^n$ in
\eqref{PFE-7}
yields that the above is smaller than $K(u_n^4+\De_n^{1-\be
/2}u_n^{\be})$,
which in turn is smaller than
$Ku_n^4$ by \eqref{P}. This yields \eqref{PFE-34}.
\end{pf}

Below, we use the simplifying notation:
%
%
\begin{equation}
\label{PFE-28}\qquad V \bigl(\ka,u,u' \bigr)^n_t=
\ua \bigl(\ka,u+u' \bigr)^n_{t}+\ua \bigl(
\ka,\bigl |u-u'\bigr | \bigr)^n_t -2\ua(
\ka,u)^n_t \ua \bigl(\ka,u'
\bigr)^n_t.
\end{equation}

%
\begin{lem}\label{L3} For all $q\geq2$, and if $u_n'\asymp u_n$, we have
%
%
\begin{eqnarray}
\label{PFE-41} %
\bigl |\E \bigl(\xi(\ka,u_n)^{1,n}_j
\mid\f_{\ka jv_n} \bigr)\bigr  |&\leq& u_n^4\rdn
\phi_n,
\nonumber\hspace*{-35pt}
\\
\biggl |\E \bigl(\xi(\ka,u_n)^{1,n}_j \xi \bigl(
\ka,u_n' \bigr)^{1,n}_j\mid
\f_{\ka jv_n} \bigr) -\frac{1}{2k_n} V \bigl(\ka,u_n,u_n'
\bigr)^n_{\ka jv_n}\biggr  | &\leq&\phi_nu_n^4
\rdn,\hspace*{-35pt}
\\
\E \bigl(\bigl |\xi(\ka,u_n)^{1,n}_j\bigr |^q
\mid\f_{\ka jv_n} \bigr)&\leq& Ku_n^4/k_n^{q/2},
\nonumber\hspace*{-35pt}
\end{eqnarray}
%
%
\begin{eqnarray}
\label{PFE-42} %
\bigl |\E \bigl(\xi(\ka,u_n)^{2,n}_j
\mid\f_{\ka jv_n} \bigr)\bigr |&\leq& u_n^4\rdn
\phi_n,
\nonumber
\\[-8pt]
\\[-8pt]
\E \bigl(\bigl |\xi(\ka,u_n)^{2,n}_j\bigr |^q
\mid\f_{\ka jv_n} \bigr)&\leq& u_n^4\rdn
\phi_n /k_n^{q/2},
\nonumber
\end{eqnarray}
%
%
\begin{eqnarray}
\label{PFE-43} %
\bigl |\E \bigl(\xi(\ka,u_n)^{3,n}_j
\mid\f_{\ka jv_n} \bigr)\bigr |&\leq& u_n^4\rdn
\phi_n,
\nonumber
\\[-8pt]
\\[-8pt]
\E \bigl(\bigl |\xi(\ka,u_n)^{3,n}_j\bigr |^q
\mid\f_{\ka jv_n} \bigr)&\leq& Kv_n \bigl(u_n^{2q}+
\De_n^{q(1-\be/2)}u_n^{q\be} \bigr).
\nonumber
\end{eqnarray}
\end{lem}

\begin{pf} In the proof, and for simplicity, we denote by $\ze(l,w)_n$
the $l$th summand in the definition of $\xi(\ka,u_n)^{w,n}_j$, for $w=1,2,3$.

Upon expanding the product $\xi(\ka,u_n)^{1,n}_j\xi(\ka,u_n')^{1,n}_j$,
\eqref{PFE-32} and \eqref{PFE-33} and successive conditioning yield\vspace*{-1pt}
\begin{eqnarray*}
&&\Biggl | \E \bigl(\xi(\ka,u_n)^{1,n}_j\xi \bigl(
\ka,u_n' \bigr)^{1,n}_j\mid
\f_{\ka jv_n} \bigr)
\\
&&\qquad{}-\frac{1}{2k_n^2}\sum_{l=0}^{k_n-1}
\E \bigl(V \bigl(\ka,u_n,u_n'
\bigr)^n_{\ka
(jk_n+l)\De_n} \mid\f_{\ka jv_n} \bigr) \Biggr |\leq
\phi_n u_n^4 \rdn.
\end{eqnarray*}
The first part of \eqref{P} and \eqref{PFE-11} also yield for $l\leq k_n$:\vspace*{-1pt}
\begin{eqnarray*}
&&\bigl |\E \bigl(V \bigl(\ka,u_n,u_n'
\bigr)^n_{\ka(jk_n+l)\De_n}- V \bigl(\ka,u_n,u_n'
\bigr)^n_{\ka jv_n}\mid\f_{\ka jv_n} \bigr)\bigr |
\\
&&\qquad\leq K
\bigl(u_n^2 v_n+u_n^{\be}
\De_n^{1-\be/2}v_n^{\be/2} \bigr) \leq
\phi_nu_n^4\rdn,
\end{eqnarray*}
the last estimate coming from \eqref{P}. We deduce the second part of
\eqref{PFE-41}. Next,
\eqref{PFE-32} and \eqref{PFE-34} yields $|\E(\ze(l,1)_n\mid\f
^n_{jk_n+l})|
\leq\phi_nu_n^4\rdn$ and\vspace*{-2pt}\break  $\E(|\ze(l,1)_n|^q\mid\f
^n_{jk_n+l})\leq Ku^4_n$,
so we have the first part of \eqref{PFE-41}, and also the last part by
the Burkholder--Gundy and H\"older inequalities.

\eqref{PFE-42} is a simple consequence of \eqref{PFE-21},
plus the Burkholder--Gundy inequality again. Finally, \eqref{PFE-11} yields
\begin{eqnarray*}
\bigl |\E \bigl(\ze(l,3)_n\mid\f_{\ka jv_n} \bigr) \bigr |&\leq& K
\bigl(u_n^2v_n+\De_n^{1-\be/2}u_n^{\be}v_n^{\be/2}
\bigr),
\\
\E \bigl(\bigl |\ze(l,3)_n\bigr |^q\mid\f_{\ka jv_n} \bigr)&
\leq& Kv_n \bigl(u_n^{2q}+\De_n^{q(1-\be/2)}u_n^{q\be}
\bigr).
\end{eqnarray*}
Then \eqref{P} yields \eqref{PFE-43}.
\end{pf}

%
\begin{lem}\label{L4} For all $q\geq2$, and if $u_n'\asymp u_n$, we have\vspace*{-1pt}
%
%
\begin{eqnarray}
\label{PFE-44} %
&&\bigl |\E \bigl(\xi(\ka,u_n)^n_j
\mid\f_{\ka jv_n} \bigr)\bigr |\leq u_n^4\rdn\phi
_n,
\nonumber
\\[-1pt]
&&\biggl |\E \bigl(\xi(\ka,u_n)^n_j \xi \bigl(
\ka,u_n' \bigr)^n_j\mid
\f_{\ka jv_n} \bigr)
\nonumber
\\[-8pt]
\\[-8pt]
&&\qquad{}- \frac{1}{2k_n} \frac{V(\ka,u_n,u_n')^n_{\ka jv_n}} {
\ua(\ka,u_n)^n_{\ka jv_n} \ua(\ka,u_n')^n_{\ka jv_n}} \biggr | \leq
u_n^4\rdn\phi_n,\nonumber
\\
&&\E \bigl(\bigl |\xi(\ka,u_n)^n_j\bigr |^q
\mid\f_{\ka jv_n} \bigr)\leq K \biggl(\frac
{u_n^4}{k_n^{q/2}}
+u_n^{2q}v_n+\De_n^{q(1-\be/2)}u_n^{q\be}v_n
\biggr).
\nonumber
\end{eqnarray}
\end{lem}

\begin{pf} In view of \eqref{PFE-701} and of the previous lemma, the
first and last parts of \eqref{PFE-44} are obvious. For the second
part, by
virtue of the second estimate in \eqref{PFE-41}, it is enough to prove that
\[
\bigl |\E \bigl(\xi(\ka,u_n)^{z,n}_j \xi \bigl(
\ka,u_n' \bigr)^{w,n}_j\mid
\f_{\ka
jv_n} \bigr)\bigr |\leq u_n^4\rdn\phi_n
\]
for all $z,w=1,2,3$ but $z=w=1$. This property follows from the Cauchy--Schwarz
inequality and all estimates in the previous lemma with $q=2$, except when
$z=1$ and $w=3$ or $z=3$ and $w=1$.

We will examine the case $z=1$ and $w=3$, the other one being analogous.
We have
\[
\xi \bigl(\ka,u_n' \bigr)^{3,n}_j=
\frac{1}{k_n} \sum_{l=0}^{k_n-2}(k_n-l-1)
\bigl(\ua \bigl(\ka,u_n' \bigr)_{\ka(jk_n+1+l)\De_n}^n-
\ua \bigl(\ka,u_n' \bigr)^n_{\ka
(jk_n+l)\De_n}
\bigr),
\]
yielding
\[
\xi(\ka,u_n)^{1,n}_j
\xi \bigl(\ka,u_n' \bigr)^{3,n}_j=
\frac{1}{k_n^2} \sum_{l=0}^{k_n-1}\sum
_{l'=0}^{k_n-2} \bigl(k_n-l'-1
\bigr)\al^n_{l,l'},
\]
where
\begin{eqnarray*}
\al^n_{l,l'}&=& \bigl(\cos \bigl(u_n\rho(
\ka)^n_{1+\ka jk_n+\ka l} \bigr) -\ua(\ka,u_n)^n_{\ka(jk_n+l')\De_n}
\bigr)
\\
&&{}\times\bigl( \ua \bigl(\ka,u_n' \bigr)^n_{\ka(jk_n+1+l')\De_n}-
\ua \bigl(\ka,u_n' \bigr)^n_{\ka
(jk_n+l')\De_n}
\bigr) 
\end{eqnarray*}
and it is thus enough to prove that
$a^n_{l,l'}=\E(\al^n_{l,l'} \mid\f_{\ka jv_n})$ satisfies
%
%
\begin{equation}
\label{PFE-47} \bigl |a^n_{l,l'}\bigr |\leq\cases{ %
u_n^4\displaystyle\frac{\rdn}{k_n} \phi_n &
\quad$\mbox{if } l\neq l'$,\vspace*{3pt}
\cr
u_n^4 \rdn
\phi_n&\quad$\mbox{if } l=l'$. 
}
\end{equation}

If $l<l'$, and since $|\ua(\ka,u_n)^n_t|\leq1$, \eqref{PFE-11} with
$s=\De_n$
and the first part of \eqref{P} give us
\[
\bigl |\E \bigl(\al^n_{l,l'}\mid\f_{\ka(jk_n+l')\De_n} \bigr)\bigr |\leq K
\De_nu_n^{\be} \bigl |\cos \bigl(u_n\rho(
\ka)^n_{1+jk_n+l} \bigr)- \ua(\ka,u_n)^n_{\ka(jk_n+l)\De_n}
\bigr |.
\]
Then \eqref{PFE-34} and the Cauchy--Schwarz inequality yield
$|a^n_{l,l'}|\leq K\De_nu_n^{2+\be}$, so \eqref{P}
again implies \eqref{PFE-47}. If $l>l'$
\eqref{PFE-32} yields
\[
\bigl |\E \bigl(\al^n_{l,l'}\mid\f_{\ka(jk_n+l)\De_n} \bigr)\bigr |\leq
\phi_n u_n^4\rdn \bigl |\ua \bigl(
\ka,u_n' \bigr)^n_{\ka(jk_n+1+l')\De_n} -\ua
\bigl( \ka,u_n' \bigr)^n_{\ka(jk_n+l')\De_n}\bigr  |,
\]
and \eqref{PFE-11} with $s=\ka\De_n$ and the Cauchy--Schwarz inequality yield
$|a^n_{l,l'}|\leq u_n^4\De_n\phi_n$, hence \eqref{PFE-47}.

For $l=l'$, upon using \eqref{PFE-701} and
the last part of \eqref{PFE-11}, plus \eqref{PFE-34} and the Cauchy--Schwarz
inequality and \eqref{P}, we see that it is enough to prove \eqref{PFE-47}
with $\al^n_{l,l}$ replaced by
\begin{eqnarray}
\al^{\prime n}_{l,l}= \bigl( \cos
\bigl(u_n\rho(\ka)^n_{1+\ka jk_n+\ka l} \bigr) -\ua(
\ka,u_n)^n_{\ka(jk_n+l)\De_n} \bigr) \ze^n_{\ka(jk_n+l)}\nonumber
\\
\eqntext{\mbox{where } \ze^n_i=U \bigl(
\ka,u_n' \bigr)^n_{(i+\ka)\De_n} -U \bigl(
\ka,u_n' \bigr)^n_{i\De_n}.}
\end{eqnarray}
The same type of argument, now based on the first part of \eqref{PFE-10}
and \eqref{PFE-23}, plus the property $|\cos(u_n\rho(\ka)^n_i)
-\cos(u_n\rho'(\ka)^n_i)|\leq|u_n\rho''(\ka)^n_i|\wedge1$, shows
that we can
even replace $\al^{\prime n}_{l,l}$ by
\[
\al^{\prime\prime n}_{l,l}= \psi^n_{\ka(jk_n+l)}
\ze^n_{\ka(jk_n+l)} \qquad \mbox{where } \psi^n_i=
\cos \bigl(u_n \rho'(\ka)^n_{1+i}
\bigr)-\ua( \ka,u_n)^n_{i\De_n}.
\]\eject

Observe that $\ze^n_i=\sum_{w=1}^4\be(w)^n_{i}$, where
\begin{eqnarray*}
\be(1)^n_i&=&\int_{i\De_n}^{(i+\ka)\De_n}b^{U(\ka,u_n')}_s
\, ds,\qquad \be(2)^n_i=H^{U(\ka,u_n')}_{i\De_n}(W_{(i+\ka)\De_n}-W_{i\De
_n}),
\\
\be(3)^n_i&=&\int_{i\De_n}^{(i+\ka)\De_n}
\bigl(H_s^{U(\ka,u_n')} -H_{i\De_n}^{U(\ka,u_n')} \bigr)
\,dW_s,
\\
\be(4)^n_i&=&\int_{i\De_n}^{(i+\ka)\De_n}H_s^{\prime U(\ka,u_n')}
\,dW_s'+\int_{i\De_n}^{(i+\ka)\De_n}
\int_{\R}\delta^{U(\ka,u_n')}(s,z) (\um-\un) (ds,dz).
\end{eqnarray*}
By \eqref{PFE-8}, we have $|\be(1)^n_i|\leq K\De_nu_n^{\prime 2}$. Combining
\eqref{PFE-4}, \eqref{PFE-3} and \eqref{PFE-10}, we easily check that
$\E(|\be(3)^n_i|\mid\f_{(i-1)\De_n})\leq Ku_n^{\prime2}\De_n$, hence
\[
w=1,3,\qquad0\leq l\leq k_n-1 \quad\Rightarrow\quad \E \bigl(\bigl |
\psi^n_i \be(w)^n_i\bigr |\mid
\f^n_i \bigr)\leq Ku_n^4\rdn
\phi_n.
\]
A parity argument (as in Step 4 of the proof of Lemma~\ref{L1}) shows that
$\E(\psi^n_i\be(2)^n_i\mid\f_{(i-1)\De_n})=0$ for all $i$.
Finally, with
$\g^W=\si(W_s\dvtx s\geq0)$, the independence between $W$ and
$(W',\um)$ implies
that $\E(\be(4)^n_i\mid\f_{(i-1)\De_n}\vee\g^W)=0$, whereas
$\psi^n_i$ is
$\f_{(i-1)\De_n}\vee\g^W$-measurable, hence
$\E(\psi^n_i\be(4)^n_i\mid\f_{(i-1)\De_n})=0$.
All these partial results give us the needed estimate for
$|\E(\al^n_{l,l}\mid\f_{\ka(jk_n+l)\De_n})|$, and the proof is complete.
\end{pf}

%
\begin{lem}\label{L8} For any square-integrable martingale $M$ and any random
variables $\ze^n_j$ such that $|\ze^n_j|\leq K$ and each $\ze^n_j$ is
$\f_{\ka jv_n}$-measurable, and for all $t>0$, we have
%
%
\begin{equation}
\label{PFE-50}\qquad \frac{v_n}{u_n^4 \rdn} \sum_{j=0}^{[t/\ka v_n]-1}
\E \bigl((M_{\ka(j+1)v_n}-M_{\ka jv_n}) \ze^n_j
\xi(\ka,u_n)^n_j \mid\f_{\ka jv_n}
\bigr) \toop0.
\end{equation}
\end{lem}

\begin{pf} It suffices to prove the result if we replace $\xi(\ka,u_n)^n_j$
above\vspace*{-1pt} by $\xi(\ka,u_n)^{w,n}_j$, for\vspace*{1pt} $w=1,2,3$, and in this case we
denote by
$R^{w,n}_t$ the normalized sum in \eqref{PFE-50}.

When $w=2,3$, we use the following argument: the properties of $\ze
^n_j$ and
the Cauchy--Schwarz inequality yield
%
%
\begin{eqnarray}
\label{PFE-51} \E \bigl(\bigl |R^{w,n}_t\bigr | \bigr)&\leq&
\frac{Kv_n}{u_n^4 \rdn}\Biggl( \sum_{j=0}^{[t/\ka v_n]-1}
\E \bigl(\bigl |\xi(\ka,u_n)^{w,n}_j\bigr |^2
\bigr)
\nonumber
\\[-8pt]
\\[-8pt]
&&\phantom{\frac{Kv_n}{u_n^4 \rdn}\Biggl(}{}\times \sum_{j=0}^{[t/\ka v_n]} \E
\bigl((M_{\ka(j+1)v_n}-M_{\ka jv_n})^2 \bigr)
\Biggr)^{1/2},\nonumber
\end{eqnarray}
and the last sum is equal to $\E((M_{\ka v_n([t/\ka v_n]-1)}-M_0)^2)$,
which is
bounded.\vspace*{-1pt} Then it is\vspace*{-1pt} enough to show that
$\frac{v_n^2}{u_n^8\De_n}\sum_{j=0}^{[t/\ka v_n]-1}
\E(|\xi(\ka,u_n)^{w,n}_j|^2)\to0$, which follows from Lemma~\ref
{L3} and
\eqref{P}.

When $w=1$, we write $\xi(\ka,u_n)^{1,n}_j=\psi^{\prime n}_j+\psi^{\prime\prime n}_j$, where\vspace*{-2pt}
$\psi^{\prime n}_j=\frac{1}{k_n}\times\break \sum_{l=0}^{k_n-1}\eta^n_{1+\ka(jk_n+l)}$ and
$\eta^n_i=\cos(u_n\rho(\ka)^n_i)-\E(\cos(u_n\rho(\ka)^n_i)\mid
\f_{(i-1)\De_n})$,
and we are left to prove that \eqref{PFE-50} holds with $\xi(\ka,u_n)^n_j$
replaced by $\psi^{\prime n}_j$ and by $\psi^{\prime\prime n}_j$. In both cases,
we denote by $R^{\prime {1,n}}_t$ and $R^{\prime\prime {1,n}}_t$ the corresponding normalized\vspace*{-2pt}
sums. For proving $R^{\prime\prime {1,n}}_t\toop0$, we proceed as above, that is,
we have
\eqref{PFE-51} with $\psi^{\prime\prime n}_j$ instead of $\xi(\ka,u_n)^{w,n}_j$, whereas
$|\psi^{\prime\prime n}_j|\leq\phi_nu_n^4\rdn$ by \eqref{PFE-32}, hence the
result holds.

For $R^{\prime {1,n}}_t$, we observe that, by successive conditioning,
\begin{eqnarray}
R^{\prime {1,n}}_t=\frac{v_n}{u_n^4k_n\rdn}\sum
_{j=0}^{[t/\ka v_n]-1}\ze^n_j \sum
_{l=0}^{k_n-1}\E \bigl(\Beta^n_{1+\ka(k_n+l)}
\mid\f_{\ka
jv_n} \bigr)\nonumber
\\
\eqntext{\mbox{where } \Beta^n_i=M^n_i
\cos \bigl(u_n\rho(\ka)^n_i \bigr),
M^n_i=M_{(i-1+\ka)\De_n}-M_{(i-1)\De_n}.}
\end{eqnarray}
As above, $\sum_{i=1}^{[t/\De_n]}\E((M^n_i)^2)\leq K$ and
$|\cos(u_n\rho(\ka)^n_i)-\cos(u_n\rho'(\ka)^n_i)|\leq\break  K(
|u_n\*\rho''(\ka)^n_i|\wedge1)$. Hence, if
\begin{eqnarray}
\BR^{\prime {1,n}}_t=\frac{v_n}{u_n^4k_n\rdn}\sum
_{j=0} ^{[t/\ka v_n]-1}\ze^n_j \sum
_{l=0}^{k_n-1} \E \bigl( \Beta^{\prime n}_{1+\ka(k_n+l)}
\mid\f_{\ka jv_n} \bigr)\nonumber
\\
\eqntext{\mbox{where } \Beta^{\prime n}_i=M^n_i
\cos \bigl(u_n\rho'(\ka)^n_i
\bigr),}
\end{eqnarray}
by \eqref{PFE-23} and the Cauchy--Schwarz
inequality and $|\ze^n_j|\leq K$, we have for all $\ep>0$ arbitrarily small:
\begin{eqnarray*}
\E \bigl(\bigl |R^{\prime {1,n}}_t-\BR^{\prime {1,n}}_t\bigr |
\bigr)&\leq& \frac{K\sqrt{t} v_n\De_n^{-\ep/2-\be/4}}{u_n^4k_n\rdn} \Biggl(\sum_{i=1}^{[t/\De_n]}
\E \bigl( \bigl(M^n_i \bigr)^2 \bigr)
\Biggr)^{1/2}
\\
&\leq&\frac{K\sqrt{t} \De_n^{1/2-\ep/2-\be/4}}{u_n^4} \to0
\end{eqnarray*}
[use \eqref{P} again].
Now, by classical arguments, it suffices to show that $\BR
^{\prime {1,n}}_t\toop0$ when
$M$ is orthogonal to $W$, or is equal to $W$ itself. In the second
case, we
clearly have $\E(\Beta^{\prime n}_i\mid\f^n_{i-1})=0$. In the first case,
we have the
same by an application of It\^o's formula. So $\BR^{\prime {1,n}}_t=0$ in all cases,
and the proof is complete.
\end{pf}

At this stage, we can prove Lemmas~\ref{L20}, \ref{L22} and~\ref{L25}.

\begin{pf*}{Proof of Lemma~\ref{L20}}
Using \eqref{PFE-44} with
$q>2$ and $u_n\to0$ and \eqref{P} and Markov inequality yields
\begin{eqnarray*}
\PP \bigl( \bigl(\Om(\ka,\theta u_n)_{n,t}
\bigr)^c \bigr)&\leq&\sum_{j=0}^{[t/\ka v_n]-1}
\PP \biggl(\bigl |\xi(\ka,\theta u_n)^n_j\bigr |>
\frac{1}2 \biggr)
\\
&\leq&2^{-q}\sum_{j=0}^{[t/\ka v_n]}
\E \bigl(\bigl |\xi(\ka,\theta u_n)^n_j\bigr |^q
\bigr)
\\
&\leq& Kt\phi_n,
\end{eqnarray*}
hence the claim because $\Om(\ka)_{n,t}$ is a finite union of sets
$\Om(\ka,\theta u_n)_{n,t}$.
\end{pf*}

\begin{pf*}{Proof of Lemma~\ref{L22}}
The claim
$\frac{1}{u_n^2} R^{\ka,n,\theta}_t
\toop0$ readily follows from \eqref{P} and from the last part of
\eqref{PFE-44} with $q=2$ and $q=3$.

For the second claim, we set
\[
\ze^n_j=\frac{v_n}{u_n^4\rdn} \biggl( \bigl(\xi(\ka,
\theta u_n)^n_j \bigr)^2 -
\frac{2}{k_n} f_{\ka,\theta u_n}(c_{\ka jv_n}) \biggr),\qquad
\ze^{\prime n}_j=\E \bigl(\ze^n_j \mid
\f_{\ka jv_n} \bigr).
\]
By a standard martingale argument, and since $\ze^n_j$ is
$\f_{\ka(j+1)v_n}$-measurable, it is enough to show that
%
%
\begin{equation}
\label{PFE-48} \sum_{j=0}^{[t/\ka v_n]-1}
\ze^{\prime n}_j\toucp0,\qquad \sum
_{j=0}^{[t/\ka v_n]-1}\E \bigl(\bigl |\ze^n_j\bigr |^2
\bigr)\toop0.
\end{equation}

Recall that $f_{\ka,\theta u}(x)\leq Ku^4$ when $|x|\leq K$, hence
\eqref{PFE-44}
yields
\[
\E \bigl(\bigl |\ze^n_j\bigr |^2 \bigr)\leq
\frac{Kv_n^2}{\De_n} \biggl(\frac{1}{k_n^2u_n^4} +v_n+v_n
\De_n^{4-2\be}u_n^{4\be-8} \biggr).
\]
The right-hand side above is easily seen to be o$(v_n)$ by \eqref{P},
hence the second part of \eqref{PFE-48}. For the first part, we use
\eqref{PFE-44} again and also \eqref{PFE-28} and $\ua(\ka,0)^n_t=1$ to
observe that it suffices to prove that
%
%
\begin{eqnarray}
\label{PFE-49}
\qquad&&\frac{v_n}{u_n^4k_n\rdn}\sum_{i=0}^{[t/\ka v_n]-1}
\biggl(\frac{\ua(\ka,2\theta u_n)^n_{\ka jv_n}+1
-2(\ua(\ka,\theta u_n)^n_{\ka jv_n})^2} {
2(\ua(\ka,\theta u_n)^n_{\ka jv_n})^2}
\nonumber
\\[-8pt]
\\[-8pt]
\qquad&&\hspace*{197pt}{}-2f_{\ka,\theta u_n}(c_{\ka jv_n}) \biggr)\toucp0.\nonumber
\end{eqnarray}
Now we recall that $|\ua(\ka,\theta u_n)^n_t-U(\ka,\theta
u_n)_t|\leq
Ku_n^{\be}\De_n^{1-\be/2}$ and\break  $1/\ua(\ka,\theta u_n)^n_t\leq K$:
since we have
$\De_n^{1-\be/2}u_n^{\be-4}/k_n\rdn\to0$ by \eqref{P}, we can thus
substitute $\ua(\ka,\theta u_n)^n$ in \eqref{PFE-49} with $U(\ka
,\theta u_n)$. But
in this case, and by definition of $f_{\ka,u}$, each summand is identically
$0$, hence \eqref{PFE-49} is proved.
\end{pf*}

\begin{pf*}{Proof of Lemma~\ref{L25}}
Set $\Te'=\Te\backslash\{1\}
$ and
\[
\al^{\ka,n,\theta}_j= \frac{2 v_n}{(\theta u_n)^2\rdn} \biggl(\frac{2}{\ka k_n(\theta u_n)^2}
f_{\ka,\theta u_n}'(c_{\ka jv_n})-1 \biggr)
\]
and
\begin{eqnarray}
\wY^{\ka,n,\theta}_t=\sum_{j=0}^{[t/\ka v_n]-1}
\ze_j^{\ka
,n,\theta},\nonumber
\\
\eqntext{\ze_j^{\ka,n,\theta}=
\cases{ %
\al^{\ka,n,1}_j\xi(
\ka,u_n)^n_j&\quad$\mbox{if } \theta=1$,\vspace*{3pt}
\cr
\displaystyle\frac{1}{u_n^2} \bigl(\al^{\ka,n,\theta}_j\xi(\ka ,
\theta u_n)^n_j -\al^{\ka,n,1}_j
\xi(\ka,u_n)^n_j \bigr)&\quad$\mbox{if }
\theta \in\Te'$. 
}}
\end{eqnarray}
The claim of the lemma is then equivalent to saying that
$(\wY^{\ka,n,\theta})_{\theta\in\Te}$ converges stably in law to
$ (\ka^{1/2}Z,(\ka^{1/3}(\theta^2-1)\BZ)_{\theta\in\Te'} )$.

We observe that the variable $\ze_j^{\ka,n,\theta}$ is
$\f_{\ka(j+1)v_n}$-measurable, whereas \eqref{PFE-44} and \eqref{P} and
Lemma~\ref{L8} imply, for all $t>0$ and all square-integrable
martingale $M$:
\begin{eqnarray*}
\sum_{j=0}^{[t/\ka v_n]-1} \E \bigl(
\ze_j^{\ka,n,\theta}\mid\f_{\ka jv_n} \bigr)&\toop&0,
\\
\sum_{j=0}^{[t/\ka v_n]-1} \E \bigl( \bigl(
\ze_j^{\ka,n,\theta} \bigr)^4\mid\f_{jv_n}
\bigr) &\toop&0,
\\
\sum_{j=0}^{[t/\ka v_n]-1}\E \bigl(
\ze_j^{\ka,n,\theta} (M_{\ka(j+1)v_n}-M_{\ka jv_n})\mid
\f_{\ka jv_n} \bigr)&\toop&0.
\end{eqnarray*}
Hence, Theorem~2.2.15 of \cite{JP} shows that it remains to prove the
following convergences:
%
%
\begin{eqnarray}
\label{PFE-71} &&\sum_{j=0}^{[t/\ka v_n]-1}\E \bigl(
\ze_j^{\ka,n,\theta} \ze_j^{\ka
,n,\theta'} \mid
\f_{\ka jv_n} \bigr)
\nonumber
\\[-8pt]
\\[-8pt]
&&\qquad\toop\Ga^{\ka,\theta,\theta'}_t=\cases{
2\ka\displaystyle\int_0^tc_s^2
\,ds&\quad$\mbox{if } \theta=\theta'=1$,\vspace*{3pt}
\cr
0&\quad$\mbox{if }
\theta=1\neq\theta'$,\vspace*{3pt}
\cr
\displaystyle\frac{\ka^3}6 \bigl(
\theta^2-1 \bigr) \bigl(\theta^{\prime2}-1 \bigr)\int
_0^tc_s^4 \,ds &
\quad$\mbox{if } \theta,\theta'\in\Te'$. 
}\nonumber
\end{eqnarray}

Recalling $|f_{\ka,\theta u_n}'(c_t)|\leq Ku_n^4$, it is enough to show
that $\Ga_t^{\ka,n,\theta,\theta'}\toop\Ga_t^{\ka,\theta,\theta
'}$, where
\[
\Ga_t^{\ka,n,\theta,\theta'}=\cases{ %
\displaystyle
\frac{4v_n^2}{u_n^4\De_n}\sum_{j=0}^{[t/\ka v_n]-1} \E \bigl(
\bigl(\xi(\ka,u_n)_j \bigr)^2\mid
\f_{\ka jv_n} \bigr)\vspace*{3pt}
\cr
\qquad\mbox{if } \theta =\theta '=1,\vspace*{3pt}
\cr
\displaystyle\frac{4v_n^2}{u_n^6\De_n}\sum_{j=0}^{[t/\ka v_n]-1} \E
\biggl(\xi(\ka,u_n)^n_j \biggl(
\frac{\xi(\ka,\theta'
u_n)^n_j}{\theta^{\prime2}} -\xi(\ka,u_n)^n_j
\biggr)\Bigm|\f_{\ka jv_n} \biggr)\vspace*{3pt}
\cr
\qquad\mbox{if } \theta=1\neq
\theta',\vspace*{3pt}
\cr
\displaystyle\frac{4v_n^2}{u_n^8\De_n}\sum_{j=0}^{[t/\ka v_n]-1}
\E \biggl( \biggl(\frac{\xi(\ka,\theta u_n)^n_j}{\theta^2} -\xi(\ka,u_n)^n_j
\biggr)\vspace*{3pt}
\cr
\hspace*{82pt}{}\times \biggl(\displaystyle\frac{\xi(\ka,\theta' u_n)^n_j}{\theta^{\prime2}} -\xi(\ka,u_n)^n_j
\biggr)\Bigm|\f_{\ka jv_n} \biggr) \vspace*{3pt}
\cr
\qquad\mbox{if } \theta,
\theta'\in\Te'. 
}
\]
We\vspace*{-1pt} then apply \eqref{PFE-44} again, plus $v_n=k_n\De_n$ and the fact that
$v_n/u_n^4\rdn\to0$ by \eqref{F-8}, and conclude that it is enough
to show
$\Ga_t^{\prime\ka,n,\theta,\theta'}\toop\Ga_t^{\ka,\theta,\theta
'}$, where
\[
\Ga_t^{\prime\ka,n,\theta,\theta'}\!=\cases{ %
\displaystyle
\frac{2v_n}{u_n^{4}}\sum_{j=0}^{[t/\ka v_n]-1}
\frac{V(\ka,u_n,u_n)^n_{\ka jv_n}}{(\ua(\ka,u_n)^n_{\ka jv_n})^2}\vspace*{3pt}
\cr
\qquad\mbox{if } \theta=\theta'=1,\vspace*{3pt}
\cr
\displaystyle\frac{2v_n}{u_n^6}\sum_{j=0}^{[t/\ka v_n]-1}
\biggl( \frac{V(\ka,u_n,\theta' u_n)^n_{\ka jv_n}}{\theta^{\prime2} \ua(\ka
,u_n)^n_{\ka jv_n}
\ua(\ka,\theta' u_n)^n_{\ka jv_n}} -\frac{V(\ka,u_n,u_n)^n_{\ka jv_n}}{(\ua(\ka,u_n)^n_{\ka
jv_n})^2} \biggr)\vspace*{3pt}
\cr
\qquad\mbox{if }
\theta=1\neq\theta',\vspace*{3pt}
\cr
\displaystyle\frac{2v_n}{u_n^8}\sum
_{j=0}^{[t/\ka v_n]-1} \biggl( \frac{V(\ka,\theta u_n,\theta' u_n)^n_{\ka jv_n}} {
\theta^2\theta^{\prime2} \ua(\ka,u_n)^n_{\ka jv_n} \ua(\ka,\theta'
u_n)^n_{\ka jv_n}} +
\frac{V(\ka,u_n,u_n)^n_{\ka jv_n}}{(\ua(\ka,u_n)^n_{\ka jv_n})^2}\vspace*{3pt}
\cr
\hspace*{59pt}{}-\displaystyle\frac{V(\ka,u_n,\theta u_n)^n_{\ka jv_n}}{\theta^2
\ua(\ka
,u_n)^n_{\ka jv_n}
\ua(\ka,\theta' u_n)^n_{\ka jv_n}}\vspace*{3pt}
\cr
\hspace*{147pt}{}-\displaystyle\frac{V(\ka,u_n,\theta' u_n)^n_{\ka jv_n}}{\theta^{\prime2} \ua(\ka
,u_n)^n_{\ka jv_n}
\ua(\ka,\theta' u_n)^n_{\ka jv_n}}
\biggr)\!\vspace*{3pt}
\cr
\qquad\mbox{if } \theta ,\theta '\in\Te'.
}
\]

If we denote $\Ga_t^{\prime\prime \ka,n,\theta,\theta'}$ the same as above, with
$\ua(\ka,u)^n$ and $V(\ka,u,u')^n$ substituted with $U(\ka,u)$ and
$U(\ka,u+u')+U(\ka,|u-u'|)-2U(\ka,u)U(\ka,u')$, and upon using
\eqref{PFE-701}
and $|\ua(\ka,u)^n_t-U(\ka,u)_t|\leq Ku^{\be}\De_n^{1-\be/2}$ and the
same argument as in the proof of the previous lemma, we see that it remains
to prove $\Ga_t^{\prime\prime \ka,n,\theta,\theta'}\toop\Ga_t^{\ka,\theta
,\theta'}$. The form
\eqref{PFE-7} of $U(\ka,u)$ allows us to check that indeed
\[
\Ga_t^{\prime\prime \ka,n,\theta,\theta'}=\cases{ %
\displaystyle
\frac{8v_n}{u_n^{4}}\sum_{j=0}^{[t/\ka v_n]-1}
f_{\ka,u_n}(c_{\ka jv_n}) \hspace*{94pt}\qquad\mbox{if } \theta=
\theta'=1,\vspace*{3pt}
\cr
\displaystyle \frac{8v_n}{u_n^6}\sum
_{j=0}^{[t/\ka v_n]-1} \biggl( \frac{f_{\ka,u_n\sqrt{\theta'}} (c_{\ka jv_n})}{\theta^{\prime2}}
-f_{\ka,u_n}(c_{\ka jv_n}) \biggr)\qquad\mbox{if } \theta=1\neq
\theta',\vspace*{3pt}
\cr
\displaystyle\frac{8v_n}{u_n^8}\sum
_{j=0}^{[t/\ka v_n]-1} \biggl( \frac{f_{\ka,u_n\sqrt{\theta\theta'}} (c_{\ka jv_n})}{\theta
^2\theta^{\prime2}}
+f_{\ka,u_n}(c_{\ka jv_n})\vspace*{3pt}
\cr
\hspace*{60pt}{}-\displaystyle\frac{f_{\ka,u_n\sqrt{\theta}} (c_{\ka
jv_n})}{\theta^2} -
\frac{f_{\ka,u_n\sqrt{\theta'}} (c_{\ka jv_n})}{\theta^{\prime2}} \biggr)
\cr
\hspace*{203pt}\qquad\mbox{if } \theta,\theta'\in
\Te'. 
}
\]

Observing that $ |f_{\ka,y}(x)-\frac{\ka^2}4 y^4x^2
-\frac{\ka^4}{48} y^8x^4 |\leq
Ky^{12}$ for all $x,y$ within\vspace*{-2pt} an arbitrary compact set, we readily obtain
$\Ga_t^{\prime\prime \ka,n,\theta,\theta'}\toop\Ga_t^{\ka,\theta,\theta'}$
from a Riemann sum approximation and $u_n\to0$. This completes the proof.
\end{pf*}

\subsection{Proof of Theorem \texorpdfstring{\protect\ref{TEF-2}}{5}}

At this stage, Theorem~\ref{TEF-2} is the only result left to be proven. In view of \eqref{F-5}
and by expanding $x\mapsto\log(\cos x)$ near $0$ and using the
boundedness of the process $a_t'$, we get the following
bound, uniform in $u\in(0,1]$:
\[
\bigl |A'(u)_t^n-A(u)^n_t\bigr |
\leq Ktu^{2\be-2} \De_n^{2-\be},
\]
implying
\[
\frac{1}{u_n^2 \rdn} \bigl(A'(\theta u_n)^n-A(
\theta u_n)^n \bigr) \toucp0
\]
if $\be<\frac{3}2$ because of \eqref{P}. Recall also that $A'(u)^n=A(u)^n$
when $\gamma^++\gamma^-=0$ identically. Henceforth, if we put
\[
\WZ(\ka,u)^n_t=\frac{1}{\rdn} \bigl(\wC(\ka
,u)^n_t-C_t-A(u)^n_t
\bigr),
\]
we have the following consequence of Theorem~\ref{TF-1}:
Under the assumptions of this theorem, then
%
%
\begin{eqnarray}
\label{EF-1}
&&\biggl(\WZ(\ka,u_n)^n, \biggl(
\frac{1}{u_n^2} \bigl(\WZ(\ka,\theta u_n) -\WZ(
\ka,u_n)^n \bigr) \biggr) _{\theta\in\Te} \biggr)
\nonumber
\\[-8pt]
\\[-8pt]
&&\qquad\tolls
\bigl(\ka^{1/2}Z, \bigl(\ka^{3/2} \bigl(\theta ^2-1
\bigr)\BZ \bigr) _{\theta\in\Te} \bigr)\nonumber
\end{eqnarray}
for $\ka=2$, and also for $\ka=1$ when either $1<\be<\frac{3}2$ or
$\be\geq\frac{3}2$ and $\gamma^++\gamma^-=0$ identically, that is,
under the
conditions of Theorem~\ref{TEF-2}.\vadjust{\goodbreak}

We choose a number $\ze>1$, and observe that $\wC(u,\ze)^n_T=
\wC(2,u,\ze)^n_T$ and $\wC'(u,\ze)^n_T=\wC(1,u,\ze)^n_T$, where
\[
\wC(\ka,u,\ze)^n_T=\wC(\ka,u)^n_T
-\frac{(\wC(\ka,\ze u)^n_T-\wC(\ka,u)^n_T)^2} {
\wC(\ka,\ze^2 u)^n_T-2\wC(\ka,\ze u)^n_T+\wC(\ka,u)^n_T}.
\]
By the definition of $A(u)^n$, we have $A(\ze u)^n_t=\ze^{\be-2}A(u)^n_t$.
Hence, with $\eta=\ze^{\be-2}-1$, we get
%
%
\begin{eqnarray}
\label{EF-0} %
 \wC(\ka,u_n,
\ze)^n_T&=&C_T+A(u_n)^n_T+
\rdn\WZ(\ka,u_n)^n_T-\frac{ (\eta A(u_n)^n_T+u_n^2\rdn\Phi_n )^2} {
\eta^2 A(u_n)^n_T+u_n^2\rdn\Phi_n'},\nonumber
\\\\[-24pt]
 \eqntext{\mbox{where } \Phi_n=\displaystyle\frac{1}{u_n^2} \bigl(\WZ(\ka,\ze
u)^n_T-\WZ(\ka,u)^n_T \bigr),}
\\
 \eqntext{\Phi_n'=\displaystyle\frac{1}{u_n^2} \bigl(\WZ \bigl(\ka,
\ze^2 u \bigr)^n_T-2\WZ(\ka,\ze
u)^n_T +\WZ(\ka,u)^n_T \bigr).}
\end{eqnarray}

Now, \eqref{EF-1} applied with $\Te=\{1,\ze,\ze^2\}$ yields
%
%
\begin{equation}
\label{EF-2} \quad\qquad\bigl(\WZ(\ka,u_n)^n_T,
\Phi_n,\Phi_n' \bigr) \tols \bigl(
\ka^{1/2}Z_T,\ka^{3/2} \bigl(\ze^2-1
\bigr)\BZ_T,\ka^{3/2} \bigl(\ze ^2-1
\bigr)^2\BZ_T \bigr).
\end{equation}
Recall also that $A(u)^n_t=u^{\be-2}\De_n^{1-\be/2}A_t$, where
$A_t=2\int_0^ta_s\,ds$. We then single out two cases:

First, on the set $\{A_T=0\}$, we have
\[
\frac{1}{\rdn} \bigl(\wC(\ka,u_n,\ze)^n_T-C_T
\bigr)=\WZ(\ka,u_n)^n_T+u_n^2
\frac{\Phi_n^2}{\Phi_n'}
\]
and \eqref{EF-2} shows that the ratio $\Phi_n^2/\Phi_n'$ converges in
law to $\ka^{3/2}\BZ_T$ ($\f$-condition\-ally Gaussian with positive variance,
hence nonvanishing almost surely). Since $u_n\to0$, another
application of
\eqref{EF-2} readily yields that, in restriction to the set $\{A_T=0\}$,
the variables $\frac{1}{\rdn}(\wC(\ka,u_n,\ze)^n_T-C_T)$ converge
stably in law
to $\ka^{1/2}Z_T$.

Second, we look at what happens on the set $\{A_T>0\}$, on which we have
by a simple calculation:
\begin{eqnarray*}
&&\frac{1}{\rdn} \bigl(\wC(\ka,u_n,\ze)^n_T-C_T
\bigr)
\\
&&\qquad=\WZ(\ka,u_n)^n_T +\frac{u_n^2(\Phi_n'-2\eta\Phi_n)A_T+u_n^{6-\be}\De_n^{(\be
-1)/2}\Phi_n^2} {
\eta^2A_T+u_n^{4-\be}\De_n^{(\be-1)/2}\Phi_n'}.
\end{eqnarray*}
Then \eqref{EF-2} again yields that, in restriction to the set $\{
A_T>0\}$,
the variables $\frac{1}{\rdn}(\wC(\ka,u_n,\ze)^n_T-C_T)$ converge
stably in law
to $\ka^{1/2}Z_T$.

So, indeed $\frac{1}{\rdn}(\wC(\ka,u_n,\ze)^n_T-C_T)\tols\ka
^{1/2}Z_T$ on
$\Om$,
which completes the proof of Theorem~\ref{TEF-2}.

\section*{Acknowledgements}

We would like to thank the Editor, an Associate Editor and two referees
for very helpful suggestions.


%

\printaddresses

\end{document}